\documentclass[11pt]{article}

\usepackage[toc,page]{appendix}
\usepackage{lipsum} 
\usepackage{url}
\usepackage{amsmath}
\usepackage{amsthm}
\usepackage{amssymb}
\usepackage{amsthm}
\usepackage{bm}
\usepackage{boxedminipage}
\usepackage{caption}
\usepackage{array}
\usepackage{graphicx}
\usepackage{float}
\usepackage{color}
\usepackage{enumerate}
\usepackage{amsthm,amssymb,amsmath}
\usepackage{geometry}
\newtheorem{theorem}{Theorem}
\newtheorem{example}{Example}
\newtheorem{lemma}{Lemma}
\newtheorem{proposition}{Proposition}

\newtheorem{definition}{Definition}

\newtheorem{remark}{Remark}
\usepackage{url}
\usepackage{amssymb}
\usepackage{bm}
\usepackage{boxedminipage}
\usepackage{array}
\usepackage{graphicx}
\usepackage{float}
\usepackage{enumerate}
\usepackage{amssymb,amsmath}
\usepackage{geometry}
\usepackage{lscape}
\usepackage{bbold}
\usepackage{multirow}
\usepackage{subcaption}
\usepackage[toc,page]{appendix}
\usepackage{endnotes}
\let\footnote=\endnote

\usepackage{ textcomp }
\newcommand{\adj}{y}
\newcommand{\nonadj}{x}
\newcommand{\adjset}{\mathcal{Y}}

\newcommand{\uncert}{\zeta}

\newcommand{\uncertset}{\mathcal{Z}}

\newcommand{\uncertA}{W}

\newcommand{\uncertb}{u}
\newcommand{\hs}{\hspace{.3cm}}
\newcommand{\trace}[2]{trace\left(#1#2\right)}
\newcommand{\SDP}[1]{S^+_#1}
\newcommand{\Eucledianorm}[1]{\left\Vert#1\right\Vert_2}
\newcommand{\scalervar}{\eta}
\newcommand{\identity}[1]{I_{#1}}
\newcommand{\zero}[1]{0_{#1}}
\newcommand{\ahmadreza}[1]{\textcolor{black}{#1}}
\title{Extending the Scope of Robust Quadratic Optimization}

\author{
	Ahmadreza Marandi\footnote{Department of Industrial Engineering and Innovation Sciences, Eindhoven University of Technology, The Netherlands} \footnote{Corresponding author: a.marandi@tue.nl}\\
Aharon Ben-Tal\footnote{CentER Extramural fellow, Tilburg University, The Netherlands}\\
Dick den Hertog\footnote{Tilburg School of Economics and Management, Tilburg University, The Netherlands}\\
Bertrand Melenberg$^\text{\textsection}$
}
\date{\today}
\begin{document}
 \captionsetup{labelfont={bf},format=hang,font=it,labelsep=period}
 
 \maketitle
 \abstract{We derive computationally tractable formulations of the robust counterparts of convex quadratic and conic quadratic constraints that are concave in matrix-valued uncertain parameters. We do this for a broad range of uncertainty sets. In particular, we show how to reformulate the support functions of uncertainty sets represented in terms of matrix norms and cones. Our results provide  extensions to known results from the literature. We also consider hard quadratic constraints; those that are convex in uncertain matrix-valued parameters. For the robust counterpart of such constraints we derive inner and outer tractable approximations. As application, we show how to construct a natural uncertainty set based on a statistical confidence set around a sample mean vector and covariance matrix and use this to provide a tractable reformulation of the robust counterpart of an uncertain portfolio optimization problem. We also apply the results of this paper to a norm approximation and a regression line problem. }
 \section{Introduction}
 Many real-life optimization problems have parameters whose values are not exactly known.  Let us consider an optimization problem containing the constraint
 \begin{equation}\label{Eq:general uncertain constraint}
 f(\adj,\zeta)\leq0, 
 \end{equation}
 where $\adj\in\mathbb{R}^n$ is the decision variable, $\zeta\in\mathbb{R}^{ t}$ is the parameter that is not known exactly, and $f:\mathbb{R}^{n}\times \mathbb{R}^{t }\longrightarrow\mathbb{R}$ is a continuous function. One way to deal with parameter uncertainty is Robust Optimization (RO), which enforces the constraints to hold for all uncertain parameter values in a user specified uncertainty set $\mathcal{Z}\subseteq\mathbb{R}^{t}$. More precisely, RO changes \eqref{Eq:general uncertain constraint} into
 \begin{equation}\label{Eq:general RO of uncertain constraint}
 f(\adj,\zeta)\leq0, \hs\forall\zeta \in\uncertset.
 \end{equation}
 This leads to a semi-infinite optimization problem, called the robust counterpart (RC), which is generally computationally intractable \ahmadreza{(see, e.g., Example 1.2.7 of the book \cite{bental2009robust})}. A challenge in RO is to find a  tractable, i.e., conic quadratic or semi-definite, reformulation of the RC. Tractability depends not only on the functions defining the constraint, i.e., $f(\adj,\zeta)$ in \eqref{Eq:general uncertain constraint}, but also on the uncertainty set $\uncertset$. For a linear constraint with linear uncertainty, where $f(\adj,\zeta)=b(\zeta)^T\adj+c$, with scalar $c\in\mathbb{R}$ and affine $b(\zeta)\in \mathbb{R}^{n}$, there is a broad range of uncertainty sets for which the RC has a tractable reformulation, see  \cite{gorissen2015practical}. 
 
 An extension of the linear case that we consider in this paper is an uncertain quadratic constraint
 \begin{equation}\label{eq:uncertain quadratic}
 \adj^T A\ahmadreza{(\Delta)}\adj +b\ahmadreza{(\Delta)}^T\adj+c \leq 0,
 \end{equation}
 where $A\ahmadreza{(\Delta)}\in\mathbb{R}^{n\times n}$ and $b\ahmadreza{(\Delta)}\in\mathbb{R}^n$ are uncertain, and $c\in\mathbb{R}$ is deterministic. We consider uncertain constraints in which the uncertainty in the parameters can be formulated in a matrix format, whereas the results in the literature are mainly for vector uncertainty.  Throughout the paper, we use the notation $\zeta$ in case of vector uncertainty and $\Delta$ in case of matrix uncertainty. So, we consider \ahmadreza{the RC of \eqref{eq:uncertain quadratic}:}
 \begin{subequations}\label{concave uncertainty}
 	\begin{eqnarray}
 	&\adj^T A(\Delta)\adj +b(\Delta)^T\adj+c \leq 0,\hs\forall \Delta \in \uncertset, \label{concave uncertainty convex quadratic}
 	\end{eqnarray}
 	where $\Delta\in\mathbb{R}^{n\times n}$ (of the same dimension as $A$) is the uncertain parameter belonging to the convex compact uncertainty set $\uncertset\subset\mathbb{R}^{n\times n}$, where $A(\Delta)\in  \mathbb{R}^{n \times n}$ and $b(\Delta)\in \mathbb{R}^n$ are affine in $\Delta$, $A(\Delta)$ is positive semi-definite for all $\Delta \in\mathcal{Z}$, and where $c \in \mathbb{R}$.
 	
 	An important optimization problem
 	having constraints in the form (\ref{concave uncertainty convex quadratic}), 
 	is a \textit{portfolio choice problem}, in which one tries to find  an asset allocation that trades off a low risk against a high expected return. One can formulate a portfolio choice problem using the form (\ref{concave uncertainty convex quadratic}), where $A(\Delta)$ is the covariance matrix and $b(\Delta)$ is minus the vector of  mean returns (possibly with a weight), respectively.
 	
 	In addition to a quadratic constraint in the form (\ref{concave uncertainty convex quadratic}), we consider a conic quadratic constraint that is concave in the uncertain parameters in the form
 	\begin{equation}
 	\sqrt{\adj^T A(\Delta)\adj} +b(\Delta)^T\adj+c \leq 0,\hs\forall \Delta \in \uncertset,\label{concave uncertainty conic quadratic}
 	\end{equation}
 	where $A(\Delta), $ $b(\Delta)$, $c$, and $\uncertset$ are defined as above.
 \end{subequations}
 
 To the best of our knowledge, there are only a few papers treating the constraints in the forms (\ref{concave uncertainty}). Moreover, the matrix $A$ typically is given as an uncertain linear combination of some primitive matrices with vector uncertainty. For example, the authors in \cite{goldfarb2003robust} study constraints in the form (\ref{eq:uncertain quadratic}), where $A$ is formulated as $\sum_{i=1}^t\zeta_iA_i$ and $\zeta=[\zeta_1,...,\zeta_t]^T \in\uncertset\subseteq \mathbb{R}^t$ is the uncertain parameter vector, for given positive semi-definite matrices $A_i$, $i=1,...,t$. They provide exact tractable reformulations of RCs for polyhedral and ellipsoidal uncertainty sets. The uncertainty set $\uncertset$ that we consider in this paper is a matrix-valued one, which is not studied in \cite{goldfarb2003robust}. The results in \cite{Jeyakumar2015} are similar to the results in \cite{goldfarb2003robust} when applied to a quadratic constraint in the form \eqref{eq:uncertain quadratic}. 
 In a more general setting, the authors in \cite{BramDual} introduce a dual problem to a general convex nonlinear robust optimization problem where the objective function and constraints are concave in the uncertain parameters, and provide conditions under which strong duality holds.

 Except for the aforementioned papers, the focus in the literature remarkably is on the constraints in the forms 
 \begin{subequations}\label{convex uncertainty}
 	\begin{eqnarray}
 	&	\adj^T A(\Delta)^TA(\Delta)\adj +b(\Delta)^T\adj+c \leq 0, \hs\forall \Delta \in \uncertset,\label{convex uncertainty convex quadratic}\\
 	&	 \sqrt{\adj^T A(\Delta)^TA(\Delta)\adj} +b(\Delta)^T\adj+c \leq 0,\hs\forall \Delta \in \uncertset,\label{convex uncertainty conic quadratic}
 	\end{eqnarray}
 \end{subequations}
 where $A(\Delta)\in\mathbb{R}^{m\times n}$ and $b(\Delta)\in\mathbb{R}^{ n}$ are affine in $\Delta\in \uncertset\subseteq\mathbb{R}^{m\times n}$, and $\uncertset$ is a convex compact set.
 For example, the book \cite{bental2009robust} and papers \cite{el1997robust} and \cite{bental2002quadratic} treat the constraints in the forms (\ref{convex uncertainty}). The drawback of (\ref{convex uncertainty}) is that the RC is, in general, (computationally) intractable, since the constraints are convex in the uncertain parameter $\Delta$ \ahmadreza{(see, e.g., \cite{Pardalos1991})}. 
 
 It is worth mentioning that the key characteristic of the constraints in the forms (\ref{concave uncertainty}) is that they are concave in $\Delta$ with convex $\uncertset$. Constraints in the forms (\ref{convex uncertainty}) can be formulated in terms of (\ref{concave uncertainty})\ahmadreza{, for instance as follows}:
 $$
 \begin{aligned}
 \adj^T B(\bar{ \Delta}) \adj +b^T\adj+c \leq 0, \hs\forall \bar{ \Delta} \in \bar{ \uncertset},\\
 \sqrt{\adj^T B(\bar{ \Delta})\adj} +b^T\adj+c \leq 0,\hs\forall \bar{ \Delta}\in \bar{\uncertset},
 \end{aligned}
 $$
 where $B(\bar{\Delta})=\ahmadreza{\bar{\Delta}}$ and $\ahmadreza{\bar{ \uncertset}=\left\{\bar{ \Delta}:\enskip \bar{ \Delta}=A(\Delta)^TA(\Delta),\;\;\Delta\in\uncertset  \right\}}$, \ahmadreza{but} $\bar{ \uncertset}$ is not convex \ahmadreza{anymore}, even \ahmadreza{not} for a convex $\uncertset$.
 
 \ahmadreza{On the one hand, the focus of the literature is on reformulating the RCs of constraints in the forms (\ref{convex uncertainty}) with specific convex compact uncertainty sets, with applications specially in least-squares problems. On the other hand, many applications that naturally contain constraints in the forms (\ref{concave uncertainty}) with matrix-valued uncertainty sets have been left out from the literature.
 	Some of the applications, in addition to portfolio choice problems, are the following ones.}
 \begin{itemize}  
 	\item\textbf{Chance Constraint}  \cite[Chapter 1]{bazaraa2013nonlinear}:	
 	Consider a normally distributed random vector $a\in \mathbb{R}^n$. Let $y \in \mathbb{R}^n$ be the vector of decision variables and $c\in \mathbb{R}$ be a constant scalar, respectively. Then, the \textit{chance constraint} $Prob(a^Ty+c\geq 0)\geq \alpha$ is equivalent to $0\geq z_\alpha \sqrt{y^T\Sigma y}-y^T\mu-c$, where $\alpha \in (0,1)$, $z_\alpha$ is the $\alpha$ percentile of the standard normal distribution, $\mu$ and $\Sigma$ are the mean vector and covariance matrix of $a$, respectively. Usually, $\mu$ and $\Sigma$ are estimated 
 	based on historical data, which results in estimation inaccuracy. Since, $\mu$ and $\Sigma$ are uncertain, the inequality is of the format \eqref{concave uncertainty conic quadratic}. 
 	
 	\item\textbf{Quadratic Approximations}: Many optimization methods, like (quasi) Newton and Sequential Quadratic Programming, use quadratic approximations of objective and constraint functions. For a twice differentiable function, this approximation can be taken using the second order truncated Taylor expansion, which requires calculating the gradient vector and the Hessian matrix. However, often the calculated gradients and Hessians are inaccurate, which make them uncertain. 	Therefore, if we apply methods, like the Newton method, to a convex optimization problem, then we could approximate it by a convex quadratic optimization problem, with an uncertain gradient vector and Hessian matrix.  
 \end{itemize}

 The contribution of this paper is \ahmadreza{fourfold}. First, \ahmadreza{we extend the results in \cite{FenchelDuality}, who consider vector uncertainty, to derive reformulations of the support functions of matrix-valued uncertainty sets. } We derive explicit formulas for support functions of many choices of $\uncertset$, mostly of those given in terms of matrix norms and cones. 
 \ahmadreza{We demonstrate that these derivations for support functions of matrix-valued uncertainty sets are also useful for a class of linear Adjustable Robust Optimization problems introduced in \cite{ben2004adjustable}.}
 
 \ahmadreza{Second,} we derive tractable  formulations of the RCs of uncertain constraints in the forms (\ref{concave uncertainty}), where $A(\Delta)$ is positive semi-definite, with a general convex compact matrix-valued uncertainty set $\uncertset$, given in terms of its support function. \ahmadreza{	In the literature only for very special uncertainty sets tractable formulations have been developed, whereas the results in this paper are for a broad range of uncertainty sets.  }

 \ahmadreza{Third, we develop} inner and outer tractable approximations of the RCs of constraints in the forms (\ref{convex uncertainty}). 
 We do this by substituting the quadratic term in the uncertain parameter with upper and lower bounds that are linear in the uncertain parameter and hence are in the forms (\ref{concave uncertainty}). \ahmadreza{These results extend the literature in two ways. First of all, inner approximations for (\ref{convex uncertainty})  have been proposed in the literature only for box or 2-norm type uncertainty sets \cite{bental2002quadratic,ben1998robust} while our approach is for a much broader range of uncertainty sets. Secondly, in this paper we also derive outer approximations. Hence, we obtain both a lower and an upper bound for the optimal value of the problem. In the literature mostly inner approximations are derived. We test these approximations on norm approximation problems as well as linear regression problems with budgeted-type uncertainty set, which could not be treated using the results in the literature. Our numerical experiments show that the obtained robust solutions outperform the nominal solutions.}       
 
 \ahmadreza{Fourth}, we show how to construct a natural uncertainty set consisting of the mean vector and the (vectorized) covariance matrix by using historical data and probabilistic confidence sets. This type of uncertainty sets is important for applications such as portfolio optimization problems. We prove for this type of sets that the support function is semi-definite representable, and provide a tractable reformulation of the robust counterpart of an uncertain portfolio optimization problem.

 The remainder of the paper is organized as follows. Section \ref{Sec: prelemenries} introduces notations and definitions that are used throughout the paper. In \ahmadreza{Section \ref{Sec: support functions}, we show how to derive computationally tractable expressions for the support functions of matrix-valued sets defined by matrix norms and cones, and several composition rules, including summations, intersections, Cartesian products of sets, convexification, linear transformations, and many more.} In Section \ref{Sec:exact formulation}, we derive an exact tractable formulation for the RC of constraints in the forms (\ref{concave uncertainty}) for a \ahmadreza{general convex compact} uncertainty sets. In Section \ref{Sec: inner approximation}, we study constraints in the forms (\ref{convex uncertainty}) with a \ahmadreza{general convex compact} uncertainty set, and provide inner and outer approximations of the RCs.  Section \ref{Sec: new statistical uncertainty set} is about constructing an uncertainty set using historical information and confidence sets. In Section \ref{Sec: numerical result}, we apply the results of this paper to a portfolio choice, a norm approximation, and a regression line problem. This paper contains \ahmadreza{four} appendices. 
 Appendix \ref{Sec: appendix proofs} contains the proofs of lemmas and propositions not presented in the main text. The second appendix contains simple illustrative examples for the results in Section \ref{Sec:exact formulation}. \ahmadreza{In Appendix \ref{APP: some sets for assumption psd for all}, we show how one can check assumptions needed in Section \ref{Sec: inner approximation} to derive the approximations. Finally, Appendix \ref{Appendix:heuristic} contains a heuristic method to find worst-case scenarios, which are used in the numerical experiments to check the quality of the solutions obtained using the inner and outer approximations proposed in Section \ref{Sec: inner approximation}. }

 \section{Preliminaries}\label{Sec: prelemenries}
 In this section, we introduce the notations and definitions we use throughout the paper. We denote by $S_n$ the set of all $n\times n$ symmetric matrices, and by $\SDP{n}$ its subset of all positive semi-definite matrices. For  $A,B\in \mathbb{R}^{n\times n}$, the notations  $A\succeq B$ and $A\succ B$ are used when $A-B \in \SDP{n}$ and  $A-B \in int(\SDP{n})$, respectively, where $int(\SDP{n})$ denotes the interior of $\SDP{n}$. We denote by $\trace{A}{}$ the trace of $A$. For  $A,B\in \mathbb{R}^{n\times m}$, we set
 $vec(A):=[A_{11},...,A_{1m},...,A_{n1},...,A_{nm}]^T,$
 and hence, $\trace{A}{B^T}=vec(A)^Tvec(B)$. For symmetric matrices $A,B\in S_{n}$, we set
 $svec(A):=[A_{11},\sqrt{2}A_{12},...,\sqrt{2}A_{1n},A_{22},...,\sqrt{2}A_{(n-1)n},A_{nn}]^T,$
 and hence, $\trace{A}{B}=svec(A)^Tsvec(B)$.
 Additionally, to represent a vector $d\in \mathbb{R}^n$ by its components, we use $[d_i]_{i=1,...,n}$.
 Also, we denote the zero matrix in $\mathbb{R}^{n\times m}$ and identity matrix in $S_n$ by $\zero{n\times m}$ and $\identity{n}$, respectively. \ahmadreza{Moreover, for matrices $A,B\in\mathbb{R}^{n\times m}$, we denote the Hadamard product by $A\circ B;$ i.e., for any $i=1,...,n$ and $j=1,...,m,$ we have $\left(A\circ B\right)_{ij}=A_{ij}B_{ij}$.}
 
 We denote the singular values of a matrix $A\in \mathbb{R}^{m\times n}$ with rank $r$ by $\sigma_1(A)\geq...\geq\sigma_r(A)>0$.
 For a vector $x\in \mathbb{R}^n$, the Euclidean norm is denoted by $\Eucledianorm{x}$. We use the following matrix norms in this paper:
 \begin{description}
 	\item[Frobenius norm:] 
 	$
 	\Vert A\Vert_F=\sqrt{\sum_{i=1}^m\sum_{j=1}^nA_{ij}^2};
 	$
 	
 	\item[$l_1$ norm:] 
 	$
 	\Vert A\Vert_1=\sum_{i=1}^m\sum_{j=1}^{n}|A_{ij}|;
 	$
 	\item[$l_\infty$ norm:] 
 	$
 	\Vert A\Vert_\infty=\max_{1\leq i\leq m \atop 1\leq j\leq n }|A_{ij}|;
 	$
 	\item[spectral norm:] 
 	$
 	\Vert A\Vert_{2,2}=\sup_{ \Eucledianorm{x}  = 1} \Eucledianorm{Ax};
 	$ 
 	\item[trace (nuclear) norm:]
 	$
 	\Vert A\Vert_\Sigma=\sigma_1(A)+...+\sigma_r(A);
 	$
 	\item[dual norm:] For a general matrix norm $\Vert .\Vert$, its dual norm is defined as
 	$
 	\Vert A\Vert^*=\max_{\Vert B\Vert=1}\trace{B^T}{A}.
 	$
 \end{description}
 \begin{remark}
 	Let $\Vert . \Vert$ be a general vector norm. Then a matrix norm can be defined as $\Vert vec(A)\Vert$ for a matrix $A \in \mathbb{R}^{m\times n}$. \textit{Frobenius}, \textit{$l_1$}, and \textit{$l_\infty$ norms} are examples of this type of matrix norms.
                                                                 	\hfill \qed
 \end{remark}
 The following lemma provides the exact formulations of the dual norms corresponding to the matrix norms defined above. 
 \begin{lemma}\cite[Section 5.6]{horn2012matrix}\label{Lemma: dual norm}
 	\\
 	(a)  $\Vert A\Vert_F^*=\Vert A\Vert_F= \Eucledianorm{vec(A)}$;\hs\hs  (b)  $\Vert A\Vert_1^*=\Vert A\Vert_\infty$;\hs\hs (c) $\Vert A\Vert^* _\Sigma=\Vert A\Vert_{2,2}=\sigma_1(A)$. 
 	\hfill	\qed

 \end{lemma}
 In the rest of this section, we recall some definitions related to optimization.
 
 \begin{definition}\label{inner approximation}
 	Let $\adjset$ be a set determined by constraints in a variable $\adj$. A set $\mathcal{S}$ determined by constraints in the variable $\adj$ and additional variable $\nonadj$, is an inner approximation of $\adjset$, if
 	$
 	(\nonadj,\adj) \in\mathcal{S}\; \Rightarrow \; \adj\in \adjset.
 	$
 	A set $\mathcal{S}$ is an outer approximation if 
 	$
 	\adj\in \adjset\; \Rightarrow \;\exists \nonadj: (\nonadj,\adj) \in \mathcal{S}.	
 	$
 	\hfill \qed
 \end{definition}
 In  \cite{bental2009robust} the inner approximation is called safe approximation. 
 \begin{definition}\label{Def: support function}	
 	For a convex set $\uncertset$, the support function $\delta^*_\uncertset(.)$ is defined as follows:
 	$$
 	\begin{aligned}
 	&\mbox{if }\uncertset\subseteq  \mathbb{R}^n,&&\delta^*_\uncertset(\uncertb):=\sup _{b\in \uncertset} \left\{\uncertb^Tb \right\},\nonumber\\
 	&\mbox{if }\uncertset\subseteq \mathbb{R}^{m\times n},&&\delta^*_\uncertset(\uncertA):=\sup _{A\in \uncertset} \left\{\trace{A}{\uncertA^T} \right\},\nonumber\\
 	&\mbox{if }\uncertset\subseteq \mathbb{R}^{m\times n}\times \mathbb{R}^n,&&\delta^*_\uncertset(\uncertA,\uncertb):=\sup _{(A,b)\in \uncertset} \left\{\trace{A}{\uncertA^T}+\uncertb^Tb \right\},\nonumber
 	\end{aligned}
 	$$	
 	where $\uncertA\in \mathbb{R}^{m\times n}$, $\uncertb\in \mathbb{R}^n$.\hfill\qed
 \end{definition}
 \begin{definition}\label{def:slater condition}
 	\ahmadreza{ Let \begin{equation}\label{Eq:uncertainty set form for slater}
 		\uncertset=\left\{\Delta\in \mathbb{R}^{m\times n}: \hs
 		\begin{matrix}\trace{C^{i^T}}{\Delta} +q^i=0,&i=1,...,I,\\ h_\ell(\Delta)\leq 0,& \ell=1,...,L,\\ g_k(\Delta)\preceq 0_{p\times p},&k=1,...,K\end{matrix} \right\},
 		\end{equation} where $p$ is a positive integer, $C^i\in\mathbb{R}^{m\times n},$ $i=1,...,I,$ $q\in\mathbb{R}^I,$ and $h_\ell:\mathbb{R}^{m\times n}\rightarrow \mathbb{R}$, $\ell=1,...,L$, and $g_k:\mathbb{R}^{m\times n}\rightarrow \mathbb{R}^{p\times p}$, $k=1,...,K,$ are  convex continuous functions. Slater condition is satisfied for $\mathcal{Z}$ if there exists $\bar{ \Delta}\in\mathcal{Z}$ such that $h_\ell(\bar{ \Delta})<0,$ for any $\ell=1,...,L,$ and $g_k(\bar{ \Delta})\prec 0_{t\times t}$, for any $k=1,...,K.$ We call $\bar{ \Delta}$ a Slater point.\hfill \qed}
 \end{definition}
 \section{Support functions for matrix-valued uncertainty sets}\label{Sec: support functions}
 
The importance of the support functions of the vector-valued sets in the area of Robust Optimization has been highlighted by \cite{FenchelDuality}, who show how to derive explicit formulas of the support functions. Support functions of the matrix-valued sets, however, are not studied in the literature despite their applicabilities in defining uncertainty in linear optimization \cite{BERTSIMAS2004510}, quadratic optimization \cite{el1997robust}, and semi-definite optimization problems \cite{el1998robust}. For instance, let us consider an  Adjustable Robust Linear Optimization (ARO) problem, introduced in \cite{ben2004adjustable}:
 	\begin{equation}\label{Eq:ARO}
 	\begin{aligned}
 	\min_{x\in\mathbb{R}^t}c^Tx+\enskip \max_{\Delta\in\uncertset}\enskip \min_{y(\Delta)\in\mathbb{R}^n}&\enskip d^Ty(\Delta)\\
 	{s.t.}&\enskip \Delta x+By(\Delta)\leq h,\\
 	&\enskip \enskip y(\Delta )\geq0,
 	\end{aligned}
 	\end{equation}
 	where $x\in\mathbb{R}^t$ is a ``here-and-now'' decision, $\Delta \in\mathbb{R}^{m\times t}$ is the uncertain parameter, $\uncertset\subseteq\mathbb{R}^{m\times t}$ is a convex compact set, $y(.)\in\mathbb{R}^n$ is a ``wait-and-see'' variable, $c\in\mathbb{R}^t$, $B\in\mathbb{R}^{m\times n},$ $d\in\mathbb{R}^n,$ and $h\in\mathbb{R}^m.$ A typical approach to approximate the ARO problem \eqref{Eq:ARO} is to restrict the ``wait-and-see'' variable $y(\Delta )$ to be affine in the uncertainty parameter. In other words, \eqref{Eq:ARO} is approximated by
 	\begin{equation}\label{Eq:SRO}
 	\begin{aligned}
 	\min_{\tau \in\mathbb{R}, \enskip x\in\mathbb{R}^t \atop {V^i\in\mathbb{R}^{m\times n} \atop u^i\in \mathbb{R}}}&\enskip c^Tx+ \tau \\
 	{s.t.}&\enskip \tau \geq \sum_{i=1}^n d_i\left(\trace{V^i}{\Delta }+u ^i\right),&\enskip \forall \Delta \in\mathcal{Z}, \\
 	&\enskip \Delta x+\left[\sum_{i=1}^nB_{ji}(\trace{V^i}{\Delta }+u^i)\right]_{j=1,...,m}\leq h,&\enskip \forall \Delta \in\mathcal{Z},\\
 	&\enskip  \left[\trace{V^i}{\Delta }+u^i\right]_{i=1,...,n}\geq0,\enskip &\enskip \forall \Delta \in\mathcal{Z},
 	\end{aligned}
 	\end{equation}
 	which is the robust counterpart of an uncertain linear optimization problem where the uncertain parameters appear in all constraints. 
 	Problem \eqref{Eq:SRO} can be solved efficiently if the support function of $\uncertset$ has a tractable reformulation.  In this section, we focus on deriving explicit tractable formulations of the support functions of matrix-valued sets.

 
 In the following lemma we provide equivalent formulations of the support functions of the sets constructed using standard composition rules.
 \begin{lemma}\label{Lemma: composition support}Let $U\in \mathbb{R}^{n \times n}$.
 	\begin{enumerate}[(i)] 
 		\item \label{vec } Let
 		$
 		\uncertset=\left\{ \Delta\in \mathbb{R}^{n \times n}:\hs vec(\Delta)\in \mathcal{U} \subset \mathbb{R}^{n^2}\right\}.
 		$
 		Then $\delta_\uncertset^*(U)=\delta^*_\mathcal{U}(vec(U))$.
 		
 		\item \label{summation} Let $\Delta^1,...,\Delta^k\in \mathbb{R}^{n\times n}$ be given. Also, let
 		$
 		\uncertset =\left\{\sum_{i=1}^k\uncert_i\Delta^i: \hspace{.2cm \uncert \in \mathcal{U}\ahmadreza{\subseteq} \mathbb{R}^{k}		}\right\}.
 		$
 		Then, 
 		$\delta_\uncertset^*(U)=\delta^*_\mathcal{U}\left(\left[\trace{\Delta^{i}}{U^T}\right]_{i=1,...,k}\right).$
 		
 		\item\label{general format of affinity} Let $L\in \mathbb{R}^{n\times t}$ and $R\in \mathbb{R}^{s\times n}$ be given, and $\uncertset=\left\{L\Delta R:\Delta\in \mathcal{U}\subseteq\mathbb{R}^{t \times s}  \right\}$. Then $\delta^*_\uncertset(U)=\delta^*_\mathcal{U} \left(L^TUR^T\right)$.
 		
 		\item\label{hadamard product} \ahmadreza{Let $L\in\mathbb{R}^{n\times n}$ be given, and $\uncertset=\left\{\Delta:\enskip L\circ \Delta\in\mathcal{U}\subseteq\mathbb{R}^{n\times n} \right\}$. Then, 
 			$$
 			\delta^*_\uncertset(U)=\left\{ \begin{matrix}
 			\delta^*_\mathcal{U}(U\circ L^{\dagger})& \mbox{ if }  U_{ij}=0 \mbox{ for any }i,j=1,...,n, \mbox{ such that } L_{ij}=0, \\
 			+\infty& \enskip \mbox{ if }  U_{ij}\neq0 \mbox{ for some }i,j=1,...,n, \mbox{ such that } L_{ij}=0 ,
 			\end{matrix} \right.
 			$$
 			where for any $i,j=1,...,n,$
 			$$L^{\dagger}_{ij}=\left\{\begin{matrix}\frac{1}{L_{ij}}&\enskip \mbox{if } L_{ij}\neq0,\\0& \enskip \mbox{ otherwise.}\end{matrix} \right.$$
 		}
 		
 		\item \label{minkovski sum} Let $\uncertset_i\subseteq \mathbb{R}^{n \times n},$ $i=1,...,k$, and let $\uncertset=\sum_{i=1}^k\uncertset_i$ be the Minkowski sum. Then $\delta^*_\uncertset(U)=\sum_{i=1}^k\delta^*_{\uncertset_i}(U)$.
 		
 		\item \label{intersection} Let  $\uncertset_i\subseteq \mathbb{R}^{n \times n},$ $i=1,...,k$, \ahmadreza{ be in the form \eqref{Eq:uncertainty set form for slater} and have a common Slater point}. Also, let $\uncertset=\bigcap_{i=1}^k\uncertset_i$. Then $\delta^*_\uncertset(U)=\min_{U^i\in \mathbb{R}^{n \times n}\atop i=1,...,k}\left\{\sum_{i=1}^k\delta^*_{\uncertset_i}(U^i):\hs \sum_{i=1}^kU^i=U\right\}$.
 		\item\label{cartesian product dependent} Let $\uncertset_i\subseteq \mathbb{R}^{n_i \times n_i},$ $U_i\in \mathbb{R}^{n_i \times n_i}$, $i=1,...,k$, and
 		$
 		\uncertset=\left\{\Delta=(\Delta_1,...,\Delta_k):\;\; \Delta_i\in \uncertset_i,i=1,...,k\right\}.
 		$
 		Then 
 		we have $\delta^*_\uncertset\left((U_1,...,U_k)\right)=\sum_{i=1}^k\delta^*_{\uncertset_i}(U_i)$.
 		\item\label{convex hull} Let $\uncertset_i\subseteq \mathbb{R}^{n \times n},$ $i=1,...,k$, be convex and $\uncertset=conv(\bigcup_{i=1}^k\uncertset_i)$ be the convex hull. Then $\delta^*_\uncertset(U)=\max_{i=1,...,k}\delta^*_{\uncertset_i}(U)$.
 	\end{enumerate}
 \end{lemma}
 \begin{proof}{Proof.} 
 	Appendix \ref{Appendix: proof of the lemma of the support function with some composition rule }.
 \end{proof}
 \ahmadreza{Lemma \ref{Lemma: composition support} shows how we can derive the support functions of sets constructed using  different composition rules without being restricted to vector-valued sets in contrast with the results in \cite{FenchelDuality}, which hold for vector-valued sets.}
 
 In the next lemma we derive \ahmadreza{explicit} tractable reformulations of the support functions of matrix-valued uncertainty sets \ahmadreza{defined by matrix norms and the cone of positive semi-definite matrices}.
 \begin{lemma}\label{Lemma: Support functions of some uncertainty sets} Let $U\in \mathbb{R}^{n \times n}$.
 	\begin{enumerate}[(a)]
 		
 		\item \label{general norm} Let
 		$
 		\uncertset=\left\{\Delta\in \mathbb{R}^{n\times n}: \hs   \Vert \Delta\Vert \leq \rho \right\},
 		$ where $\Vert. \Vert$ is a general matrix norm. Then $\delta^*_\uncertset(U)=\rho\Vert U\Vert^*$.

 		\item \label{PSD bounded set} Let 
 		$
 		\uncertset=\{\Delta:\;\; \Delta^l\preceq \Delta \preceq \Delta^u\},
 		$
 		where $\Delta^l, \ \Delta^u\in S_n $ are given such that $ \Delta^u-\Delta^l\succ \zero{n\times n}$.
 		Then 
 		$$
 		\delta_\uncertset ^*(U)=\min_{\Lambda_1,\Lambda_2} \left\{\trace{\Delta^u}{\Lambda_2}-\trace{\Delta^l}{\Lambda_1}:\;\Lambda_2-\Lambda_1=\frac{U+U^T}{2},\; \Lambda_1,\Lambda_2\succeq\zero{n\times n}\right\}.
 		$$ 
 		
 	\end{enumerate}
 \end{lemma}
 \begin{proof}{Proof.}
 	(\ref{general norm}) \ahmadreza{This} follows directly from the definition of the dual norm.\\
 	(\ref{PSD bounded set}) \ahmadreza{This} follows directly from conic duality (see Appendix \ref{Appendix proof of lemma support function of natural sets}).
 \end{proof} 
 
 Special cases of the uncertainty sets studied in Lemma \ref{Lemma: Support functions of some uncertainty sets} have been considered in the literature. The uncertainty set constructed using the \textit{Frobenius norm} is considered in \cite{el1997robust} for the constraints in the form (\ref{convex uncertainty conic quadratic}). Also, the authors of \cite{shawe2003estimating} construct an uncertainty set for the covariance matrix using the \textit{Frobenius norm}. The constraints in the forms (\ref{convex uncertainty}) with uncertainty set defined by the \textit{spectral norm} is treated in Chapter 6 of  \cite{bental2009robust}. Furthermore, the uncertainty set that we considered in Lemma \ref{Lemma: Support functions of some uncertainty sets}(\ref{PSD bounded set}) is constructed in \cite{rife2013effect} for covariance matrices. Besides, the authors of \cite{delage2010distributionally} construct an uncertainty set for the mean vector and covariance matrix, which can be formulated as an intersection of two sets that are \ahmadreza{considered} in Lemma \ref{Lemma: Support functions of some uncertainty sets}(\ref{PSD bounded set}).
 
 It is known that the \textit{$l_1$} and $l_\infty$\textit{ norms} are linear representable  and the \textit{Frobenius norm} is conic quadratic representable. The following lemma shows that the \textit{spectral} and \textit{trace norms} are semi-definite representable.
 \begin{lemma}\label{Lemma: SDP matrix norm}  Let $U\in \mathbb{R}^{n\times n}$ and $\rho\geq0$.
 	\begin{itemize}
 		\item[(i)] $\Vert U\Vert_\Sigma\leq \rho$ if and only if there exist matrices $Y\in \mathbb{R}^{n\times n}$ and $Z\in \mathbb{R}^{n\times n}$ such that
 		$$
 		\left[\begin{matrix}
 		Y&U\\U^T&Z
 		\end{matrix}\right]\succeq\zero{2n\times 2n}, \hs \trace{Y}{}+\trace{Z}{}\leq 2\rho.
 		$$
 		\item[(ii)]$\Vert U\Vert_{2,2}\leq \rho$ if and only if $\left[\begin{matrix}
 		\rho^2\identity{n}&U\\ U^T&\identity{n}
 		\end{matrix}\right]\succeq\zero{2n \times 2n}.$
 	\end{itemize}	
 \end{lemma}
 \begin{proof}{Proof.}
 	(i) See, e.g., Lemma 1 in \cite{fazel2001rank}.\\
 	(ii)
 	See, e.g., Example 8 in \cite{FenchelDuality}, or Appendix \ref{Appendix lemma SDP refomulation of trace norm}.	
 \end{proof}
 \ahmadreza{Hitherto, we have shown how to derive tractable reformulations of the support functions of matrix-valued uncertainty sets. In the next section, we show how such reformulations can be used in tractably reformulating the RC of an uncertain quadratic constraint in the forms \eqref{concave uncertainty}. }
 \section{\ahmadreza{Tractable reformulation of Robust Quadratic Optimization problems that are concave in the uncertain parameters}} \label{Sec:exact formulation}
 \ahmadreza{In this section, we assume that $A(\Delta)=A+\Delta$. We emphasize that this assumption can be made without loss of generality for $A(\Delta)$ that is affine on $\Delta$, because of Lemma \ref{Lemma: composition support}.\eqref{general format of affinity}. Moreover, from here on in the paper, we assume that the uncertainty set $\mathcal{Z}$ is defined as in \eqref{Eq:uncertainty set form for slater} and satisfies the Slater condition.}
 The \ahmadreza{next} theorem\ahmadreza{, which is the main theorem in this section,} provides reformulations of the RCs of constraints in the forms (\ref{concave uncertainty}), \ahmadreza{and asserts that their tractabilities only depend on the uncertainty sets.}
 \begin{theorem}\label{Th: RC for conic and convex quadratic}
 	Let $\uncertset\subset \mathbb{R}^{n\times n}$ be a convex, compact set.
 	Also, let $\bar{ A}\in \mathbb{R}^{n\times n}$, $a,\bar{b}\in \mathbb{R}^n$, and $c\in \mathbb{R}$ be given. For any $\Delta \in \uncertset\subseteq\mathbb{R}^{n\times n}$, let $A(\Delta)=\bar{A}+\Delta $, $b(\Delta)=\bar{b}+\Delta a$ ($m=n$ in \eqref{concave uncertainty}). Assume that $A(\Delta)$ is positive semi-definite (PSD), for all $\Delta \in \uncertset$, and that \ahmadreza{for a Slater point $\bar{ \Delta}$}, $A(\bar{ \Delta})$ is positive definite. Then:
 	\begin{enumerate}[(I)]
 		\item\label{Th: item, convex quadratic} $y\in \mathbb{R}^n$ satisfies (\ref{concave uncertainty convex quadratic}) if and only if there exists $\uncertA \in \mathbb{R}^{n \times n}$ satisfying the convex system
 		\begin{equation}\label{SDR of convex quadratic}
 		\begin{matrix}
 		\trace{\bar{A}}{\uncertA}+\bar{b}^T\adj+ c+\delta^*_{\uncertset}( \uncertA +  \adj a^T)\leq 0,\hs\hs
 		\left[\begin{matrix}
 		\uncertA& \adj\\
 		\adj^T& 1
 		\end{matrix}\right]\succeq\zero{n+1\times n+1}.\\
 		\end{matrix}
 		\end{equation}
 		\item \label{Th: item, conic quadratic}$y\in \mathbb{R}^n$ satisfies (\ref{concave uncertainty conic quadratic}) if and only if there exist $\uncertA \in \mathbb{R}^{n \times n}$ and $\scalervar \in \mathbb{R}$ satisfying the convex system
 		\begin{equation}\label{SDR of conic quadratic}
 		\begin{matrix}
 		\trace{\bar{A}}{\uncertA}+\bar{b}^T\adj+ c+\delta^*_{\uncertset}(\uncertA +  \adj a^T)+\frac{\scalervar}{4}\leq 0,\hs\hs
 		\left[\begin{matrix} \uncertA&\adj\\\adj^T&\scalervar
 		\end{matrix}\right]\succeq \zero{n+1\times n+1} .
 		\end{matrix}
 		\end{equation}
 	\end{enumerate}
 \end{theorem}
 \begin{proof}{Proof.} 
 	To prove this theorem we use the same line of reasoning as in Theorem 2 in \cite{FenchelDuality}. For any $\Delta\in \uncertset$, it is clear that $\frac{A(\Delta)+A(\Delta)^T}{2}\succeq \zero{n\times n}$ due to positive semi-definiteness of $A(\Delta)$. Also, $\adj^TA(\Delta)\adj=\adj^T\frac{A(\Delta)+A(\Delta)^T}{2}\adj$ for any $\adj \in \mathbb{R}^n$, and $\Delta\in \uncertset$. We \ahmadreza{replace} $\adj^TA(\Delta)\adj$ by $\adj^T\frac{A(\Delta)+A(\Delta)^T}{2}\adj$ in constraints (\ref{concave uncertainty}).   \\
 	(\ref{Th: item, convex quadratic})  Let $\mathcal{U}=\left\{\left(\frac{A(\Delta)+A(\Delta)^T}{2},b(\Delta)\right):\;\Delta\in \uncertset\right\}$. It is clear that $\adj \in \mathbb{R}^n$ satisfies (\ref{concave uncertainty convex quadratic}) if and only if 
 	$
 	F(\adj):=\max_{(B,b)\in \mathcal{U}}\left\{\adj^TB\adj+b^T\adj+c\right\}\leq0.
 	$
 	Setting
 	$$
 	\delta_\mathcal{U}(B,b)=\left\{\begin{matrix}
 	0&\mbox{if }(B,b)\in \mathcal{U},\\
 	+\infty&\mbox{otherwise},
 	\end{matrix} \right.
 	$$
 	we have 
 	$
 	F(\adj)=\max_{B\succeq0_{n\times n}\atop b\in \mathbb{R}^n}\left\{\adj^TB\adj+b^T\adj+c-\delta_{\mathcal{U}}(B,b)\right\}.
 	$
 	Since $B\succeq0_{n\times n}$ for all $B\in \mathcal{U}$, and \ahmadreza{$\frac{A(\bar{ \Delta})+A(\bar{ \Delta})^T}{2}$} is positive definite \ahmadreza{and lies} in the relative interior of $\mathcal{U}$, \ahmadreza{specialization of  Theorem 4.4.3 in \cite{Borwein2005Variational} to $\mathbb{R}^{n\times n}$ implies that} $F(\adj)\leq0$ is equivalent to the existence of $\uncertA \in \mathbb{R}^{n \times n}$ and $\uncertb\in \mathbb{R}^n,$ such that
 	\begin{equation}\label{Eq: RC convex quadratic with fenchel}
 	\delta_{\mathcal{U}}^*(\uncertA,\uncertb)-\inf _{\begin{tiny}\begin{matrix}
 		A \succeq \zero{n\times n} \\b \in \mathbb{R}^n
 		\end{matrix}\end{tiny}} \left\{\trace{A}{\uncertA^T}+\uncertb^Tb-\left(\adj^T A\adj+b^T\adj+c\right) \right\}\leq 0,
 	\end{equation}
 	where $\delta_{\mathcal{U}}^*(.)$ is the support function of the set $\mathcal{U}$. 
 	It follows from Definition \ref{Def: support function} that
 	\begin{eqnarray}
 	\delta^*_\mathcal{U}(\uncertA,\uncertb)&=&\sup_{\left(B,b\right)\in \mathcal{U}}\left\{\trace{B}{\uncertA^T}+\uncertb^Tb\right\}
 	=\sup_{\Delta\in \uncertset}\left\{\trace{\frac{A(\Delta)+A(\Delta)^T}{2}}{\uncertA^T}\enskip+\enskip\uncertb^Tb(\Delta)\right\}\nonumber \\
 	&=&\sup_{\Delta \in \uncertset}\left\{\trace{\left(\bar{A}+\Delta \right)}{\left(\frac{\uncertA+\uncertA^T}{2}\right)}\enskip+\enskip\uncertb^T(\bar{b}+\Delta a)\right\}\nonumber \\
 	&=&\trace{\bar{A}}{\left(\frac{\uncertA+\uncertA^T}{2}\right)}\enskip+\enskip\uncertb^T\bar{b}\enskip+\enskip\sup_{\Delta \in \uncertset}\left\{\trace{\Delta }{\left(\frac{\uncertA+\uncertA^T}{2}\right)}+\uncertb^T\Delta a\right\}\nonumber \\
 	&=&\trace{\bar{A}}{\left(\frac{\uncertA+\uncertA^T}{2}\right)}\enskip+\enskip\uncertb^T\bar{b}\enskip+\enskip\sup_{\Delta \in \uncertset}\left\{\trace{\Delta}{\left(\left(\frac{\uncertA+\uncertA^T}{2}\right) +a\uncertb^T \right)}\right\}\nonumber \\
 	&=&\trace{\bar{A}}{\left(\frac{\uncertA+\uncertA^T}{2}\right)}\enskip+\enskip\uncertb^T\bar{b}\enskip+\enskip\delta^*_\uncertset\left( \left(\frac{\uncertA+\uncertA^T}{2}\right) +  \uncertb a^T\right).\label{Eq: support of theta}
 	\end{eqnarray}
 	Also, we have
 	\begin{eqnarray}
 	&&\inf _{\begin{tiny}\begin{matrix}
 		A\succeq \zero{n\times n} \\b \in \mathbb{R}^n
 		\end{matrix}\end{tiny}} \left\{\trace{A}{\uncertA^T}+\uncertb^Tb-\left(\adj^T A\adj+b^T\adj+c\right) \right\}\nonumber\\
 	&=&\inf _{\begin{tiny}\begin{matrix}
 		A\succeq \zero{n\times n} \\b \in \mathbb{R}^n
 		\end{matrix}\end{tiny}} \left\{\trace{A}{\uncertA}+\uncertb^Tb-\left(\adj^T A\adj+b^T\adj+c\right) \right\}\nonumber\\
 	&=&-c+\inf _{\begin{tiny}
 		\begin{matrix}	A\succeq \zero{n\times n} \\ b\in \mathbb{R}^n
 		\end{matrix}\end{tiny}}\left\{\trace{A}{ \left( \uncertA-\adj \adj^T\right) }+ b^T(\uncertb-\adj) \right\}
 	= \left\{ \begin{matrix}
 	-c &\hspace{.2cm} \uncertA -\adj\adj^T \succeq \zero{n\times n} ,\ \uncertb=\adj,\\
 	-\infty & \mbox{otherwise.}
 	\end{matrix} \right. \label{Eq: concave conjugate of convex quadratic}
 	\end{eqnarray}
 	So, the fact that $\uncertA\succeq \zero{n\times n}$ implies $\frac{\uncertA+\uncertA^T}{2}=\uncertA$, and \ahmadreza{the} Schur Complement Lemma (see, e.g., Appendix A.5.5 in  \cite{boyd2004convex}), (\ref{Eq: support of theta}), and (\ref{Eq: concave conjugate of convex quadratic}) result in (\ref{SDR of convex quadratic}).
 	\\
 	(\ref{Th: item, conic quadratic}) Similar to the proof of part (\ref{Th: item, convex quadratic}) we have $y\in \mathbb{R}^n$ satisfies (\ref{concave uncertainty conic quadratic}) if and only if there exists $\uncertA \in \mathbb{R}^{n \times n}$ such that
 	\begin{equation}\label{Eq: RC conic quadratic with fenchel}
 	\delta_{\mathcal{U}}^*(\uncertA,\uncertb)-\inf _{\begin{tiny}\begin{matrix}
 		A\succeq \zero{n\times n} \\b \in \mathbb{R}^n
 		\end{matrix}\end{tiny}} \left\{\trace{A}{\uncertA^T}+\uncertb^Tb-\left(\sqrt{\adj^T A\adj}+b^T\adj+c\right) \right\}\leq 0.
 	\end{equation}
 	Analogous to the result in Section 3.4 in \cite{BramDual}, 
 	$$
 	\inf _{\begin{tiny}\begin{matrix}
 		A\succeq \zero{n\times n} \\b \in \mathbb{R}^n
 		\end{matrix}\end{tiny}} \left\{\trace{A}{\uncertA^T}+\uncertb^Tb-\left(\sqrt{\adj^T A\adj}+b^T\adj+c\right) \right\}
 	=-c-\inf _{\scalervar}\left\{\frac{\scalervar}{4}:\; \uncertb=\adj,\;\left[\begin{matrix} \uncertA&\adj\\\adj^T&\scalervar
 	\end{matrix}\right]\succeq \zero{n+1\times n+1}\right\}. 
 	$$
 	So, (\ref{Eq: RC conic quadratic with fenchel}) is equivalent to 
 	\begin{equation}\label{Eq: RC of conic quadratic with inf}
 	\delta_{\mathcal{U}}^*(\uncertA,\uncertb)+c+\inf _{\scalervar\in \mathbb{R}}\left\{\frac{\scalervar}{4}:\;\left[\begin{matrix} \uncertA&\adj\\\adj^T&\scalervar
 	\end{matrix}\right]\succeq \zero{n+1\times n+1}\right\}\leq0.
 	\end{equation}
 	In (\ref{Eq: RC of conic quadratic with inf}), the infimum is taken over a closed lower bounded set, since $\scalervar\geq0$. Hence, $\uncertA \in \mathbb{R}^{n \times n}$ and $\adj \in \mathbb{R}^{n}$ satisfies (\ref{Eq: RC of conic quadratic with inf}) if and only if there exists $\scalervar\in \mathbb{R}$ such that
 	$$
 	\begin{aligned}
 	&\trace{\uncertA}{\bar{ A}^T}+\bar{b}^T\adj+\delta^*_\uncertset( \uncertA +  \adj a^T)+c+\frac{\scalervar}{4}\leq0,\hs
 	&\left[\begin{matrix} \uncertA&\adj\\\adj^T&\scalervar
 	\end{matrix}\right]\succeq \zero{n+1\times n+1},
 	\end{aligned}
 	$$
 	which completes the proof. 
 \end{proof}
 
 
 \ahmadreza{One of the assumptions in Theorem \eqref{Th: RC for conic and convex quadratic}} is that $A(\Delta)$ is positive semi-definite for all $\Delta\in\mathcal{Z}$. This assumption is needed to guarantee convexity of the constraint. Even though checking this assumption for a general uncertainty set is intractable (Section 8.2 in \cite{bental2009robust}), there are cases for which this assumption holds. An example is when $A(\Delta)$ is a covariance matrix, which is estimated, e.g., based on historical data. Another example is when $A(\Delta)$ is the Laplacian matrix of a weighted graph, where the weights are uncertain. In these examples,  $A(\Delta)$ by construction is positive semi-definite for all possible values of the uncertain parameter $\Delta$. Besides the aforementioned examples, it is clear that if $\bar{ A}$ is positive semi-definite and $\mathcal{Z}\subseteq\SDP{n}$, then $A(\Delta)$ is positive semi-definite for all $\Delta \in\mathcal{Z}$.

 Next to the cases mentioned above, Theorem 8.2.3 in \cite{bental2009robust} provides a tractable method to check this assumption for a specific class of uncertainty sets. In the following lemma, we mention a simplified version of this theorem.
 \begin{lemma}(Theorem 8.2.3 in \cite{bental2009robust})\label{Lemma: psd assumption holds}
 	Let $\uncertset=\left\{\Delta\ : \ \Vert \Delta\Vert_{2,2}\leq \rho  \right\}\subset \mathbb{R}^{n\times n}$. Then, for a given $\bar{ A}\in \mathbb{R}^{n\times n}$, we have that  $\bar{A}+\Delta\succeq0_{n\times n}$ for any $\Delta \in \uncertset$ if and only if $\bar{ A}-\rho I_n\succeq 0_{n\times n}.$\hfill \qed
 \end{lemma}	
 

 \subsection*{Illustrative examples }
 \ahmadreza{In the rest of this section,} we derive tractable reformulations of RCs for some natural uncertain convex quadratic and conic quadratic constraints. For brevity of exposition, we provide  in Appendix \ref{App: illustrative examples} the tractable reformulation of an uncertain convex quadratic constraint where the uncertainty set is defined by the Frobenius norm, as well as an uncertain conic quadratic constraint where the uncertainty set is similar to the one proposed in \cite{delage2010distributionally}.

 The following example is for constraints in the form \eqref{eq:uncertain quadratic}  with vector uncertainty.
 \begin{example}\label{Ex: convex quadratic with linear vector uncertainty}
 	Consider
 	\begin{equation}\label{Eq: example with vector uncertainty}
 	y^T A(\zeta)y+b(\zeta)^Ty+c\leq 0, \hs\forall \zeta\in \uncertset,
 	\end{equation}
 	where 
 	$
 	A(\zeta)=\bar{A}+\sum_{i=1}^t\zeta_iA^i, \hspace{0.3cm}b(\zeta)=\bar{b}+\sum_{i=1}^t\zeta_ib^i,
 	$
 	$(A^i,b^i)\in \mathbb{R}^{n\times n}\times \mathbb{R}^n$ is given, $i=1,...,t$. This constraint is considered in \cite{goldfarb2003robust}, where $\bar{ A}$ and $A^i$, $i=1,...,t$, are positive semi-definite and where $\uncertset=\uncertset_1\times\uncertset_2$, for some polyhedral or ellipsoidal sets $\uncertset_1\subseteq\mathbb{R}^{m}$ and $\uncertset_2\subseteq\mathbb{R}^{t-m}$, $m\in\mathbb{N}$, with $A^i=0_{n\times n}$, $i=1,...,m,$ and $b^i=0_{n\times1}$, $i=m+1,...,t$ (uncertainty in $A$ is independent of the uncertainty in $b$). In this example we show how using the results of Section \ref{Sec: support functions} can extend the results of \cite{goldfarb2003robust} for general uncertainty sets, where $A(\uncert)$ is positive semi-definite for all $\uncert \in \uncertset$. \ahmadreza{Let 
 		$$
 		\uncertset=\left\{\uncert\in \mathbb{R}^{t}: \hs
 		\begin{matrix} C\uncert +q=0_{I\times 1},&\\ h_\ell(\uncert)\leq 0,& \ell=1,...,L,\\ g_k(\uncert)\preceq 0_{p\times p},&k=1,...,K\end{matrix} \right\},
 		$$ where $p$ is a positive integer, $C\in\mathbb{R}^{I\times t},$  $q\in\mathbb{R}^I$, and where $h_\ell:\mathbb{R}^{t}\rightarrow \mathbb{R}$, $\ell=1,...,L$, and $g_k:\mathbb{R}^{t}\rightarrow \mathbb{R}^{p\times p}$, $k=1,...,K,$ are  convex continuous functions.}
 	First, we assume that $\uncertset \ahmadreza{\subseteq} \mathbb{R}^t_+$, where $\mathbb{R}^t_+$ denotes the nonnegative orthant of $\mathbb{R}^t$. Also, we assume that $\bar{ A}$ and $A_i$, $i=1,...,t$, are positive semi-definite, \ahmadreza{and there is a Slater point in $\mathcal{Z}$ for which $A(\uncert)$ is positive definite}. In this case, $\adj \in \mathbb{R}^n$ satisfies (\ref{Eq: example with vector uncertainty}) if and only if there exists $v\in \mathbb{R}^t$ such that
 	\begin{equation}\label{Eq: lemmavector uncertainty for quadratic without LMI}
 	y^T \bar{ A}y+\delta^*_{ \uncertset}( v)+\bar{ b}^Ty+c\leq 0, \;v\geq \left[\adj^TA^i\adj+b^{i^T}\adj\right]_{i=1,...,t},
 	\end{equation}
 	\ahmadreza{whose proof} can be found in Appendix \ref{Appendix proof of example with nonegative uncertainty and show no LMI}. It is clear that in this case $A(\uncert)$ is positive semi-definite for all $\uncert \in \uncertset$. Moreover, as mentioned in Remark 2 of \cite{goldfarb2003robust}, if $b^i=0$, $i=1,...,t$, then, for a general uncertainty set $\uncertset$, $\adj\in \mathbb{R}^n$ satisfies (\ref{Eq: example with vector uncertainty}) if and only if 
 	$$
 	y^T\left(\bar{A}+\sum_{i=1}^t\zeta_iA_i\right)y+\bar{b}^Ty+c\leq 0,\hs\forall \uncert \in \bar{\uncertset},
 	$$
 	where $\bar{\uncertset}=\uncertset \cap \{\zeta: \zeta\geq\zero{t\times 1}\}$ and $\bar{\uncertset} \neq \emptyset$. If the uncertainty set $\uncertset$ is a polyhedron, then a tractable RC is provided in \cite{goldfarb2003robust}. For other types of uncertainty sets, like ellipsoidal uncertainty sets, deriving tractable RCs is achievable using the results in Sections \ref{Sec: support functions} and \ref{Sec:exact formulation}. Let $\bar{\uncertset} \neq \emptyset$. Then $\adj\in \mathbb{R}^n$ satisfies (\ref{Eq: example with vector uncertainty}) if and only if 
 	$$
 	y^T \bar{ A}y+\delta^*_{\bar{ \uncertset}}( v)+\bar{ b}^Ty+c\leq 0, \;v\geq \left[\adj^TA^i\adj\right]_{i=1,...,t}.
 	$$
 	This is an extension of the results of \cite{goldfarb2003robust}, since there is a broad range of uncertainty sets for which the support functions have tractable reformulations. 
 	
 	Now, for a general case where the uncertainty in $A$ and $b$ can be dependent, if $A(\uncert)$ is positive semi-definite for all $\uncert \in \uncertset$, and positive definite for a \ahmadreza{Slater point}, then by Theorem \ref{Th: RC for conic and convex quadratic}(\ref{Th: item, convex quadratic}), $\adj$ satisfies (\ref{Eq: example with vector uncertainty}) if and only if there exists $\uncertA \in \mathbb{R}^{n\times n}$ such that
 	$$
 	\begin{matrix}
 	\trace{\bar{A}}{\uncertA}+\bar{b}^T\adj+\delta_{\uncertset}^*\left(\left[\trace{ A^i}{\uncertA}+ b^{i^T}\adj\right]_{i=1,...,t}\right)+c\leq 0,&\hs
 	\left[\begin{matrix}
 	\uncertA& \adj\\
 	\adj^T& 1
 	\end{matrix}\right]\succeq \zero{n+1\times n+1}.
 	\end{matrix}\hspace{1.5cm}\qed
 	$$
 	
 \end{example}
 
 In Section \ref{Sec: new statistical uncertainty set}, we derive a natural uncertainty set for a vector that consists of the mean vector and the vectorized covariance matrix. This type of uncertainty set can be used in different applications, such as portfolio choice problems. 
 In the following example we derive a tractable reformulation of a quadratic constraint with an uncertainty set similar to the one constructed in Section \ref{Sec: new statistical uncertainty set}.
 \begin{example}\label{Ex:general statical uncertainty set}
 	Consider the uncertain quadratic constraint 
 	$$
 	\adj^T\Delta\adj+\zeta^T\adj+c\leq 0, \hs \forall\left(
 	\zeta\atop
 	svec(\Delta)
 	\right)\in \uncertset
 	$$
 	where $\uncertset=\uncertset_1\cap\uncertset_2$, and
 	$$
 	\begin{aligned}
 	&	 \uncertset_1=\left\{\left(
 	\zeta\atop
 	svec(\Delta)
 	\right)=B\nu:\; 
 	\Eucledianorm{\nu} \leq \rho, \; \nu\in \mathbb{R}^{\frac{n^2+3n}{2}} \right\},\hs
 	\uncertset_2=\left\{\left(
 	\zeta\atop
 	svec(\Delta)
 	\right): \; \zeta\in \mathbb{R}^n,\; \Delta\in \SDP{n}\right\},
 	\end{aligned}
 	$$	
 	for some invertible $B\in \mathbb{R}^{\frac{n^2+3n}{2}\times \frac{n^2+3n}{2}}$, $\rho>0$. For a fixed $\uncertA\in S_n$, by Lemma \ref{Lemma: composition support}(\ref{intersection}), and Example 4 in \cite{FenchelDuality}, 
 	$$\delta^*_\uncertset\left(u\atop svec(W)\right)=\left\{
 	\begin{aligned}
 	\min_{u^1,u^2\atop W^1,W^2}&\hs\rho\Eucledianorm{B^T\left(u^1\atop svec(W^1)\right)}+\delta^*_{\uncertset_2}\left(u^2\atop svec(W^2)\right)\\
 	\mbox{s.t.}&\hs u^1+u^2=u,\;\hs W^1+W^2=W,\;\hs W^1,W^2\in S_n.
 	\end{aligned} \right.
 	$$ 
 	Similar to the proofs of Lemmas \ref{Lemma: composition support}(\ref{vec }) and \ref{Lemma: composition support}(\ref{cartesian product dependent}), we have
 	\begin{equation}\label{Eq: generalized statical uncertainty set}
 	\delta^*_\uncertset\left(u\atop svec(W)\right)=\min_{W^1}\left\{\rho\Eucledianorm{ B^T\left(u\atop svec(W^1)\right)} :\;W^1\succeq W\right \}.
 	\end{equation}
 	It is easy to show that there exists \ahmadreza{a Slater point in} $\uncertset$. Hence, $\adj\in \mathbb{R}^n$ satisfies (\ref{Ex:general statical uncertainty set}) if and only if there exists $\uncertA\in \SDP{n}$ that satisfies
 	$
 	\begin{aligned}
 	&\rho\Eucledianorm{B^T\left(u\atop svec\left(\uncertA\right)\right)}+c\leq0,
 	&\left[\begin{matrix}
 	\uncertA&\adj\\\adj^T&1
 	\end{matrix}\right]\succeq\zero{n+1\times n+1}.
 	\end{aligned} 
 	$
 	\hfill\qed
 \end{example}
 \section{\ahmadreza{Tractable inner and outer approximations of Robust Quadratic Optimization problems that are convex in the uncertain parameters}}\label{Sec: inner approximation}
 In this section, we provide inner and outer approximations of the RCs of constraints in the forms (\ref{convex uncertainty}) by replacing the quadratic term in the uncertain parameter with suitable upper and lower bounds.
 We assume that $A(\Delta)=\bar{A}+\Delta $, $b(\Delta)=\bar{b}+D\Delta a$, for given $D\in \mathbb{R}^{n\times m}$, \ahmadreza{full column-rank matrix} $\bar{ A}\in \mathbb{R}^{m\times n}$, \ahmadreza{vectors} $\bar{ b}\in \mathbb{R}^n$, $a\in \mathbb{R}^n$, \ahmadreza{and scalar} $c\in \mathbb{R}$, and \ahmadreza{that} $\Delta\in\uncertset\subseteq \mathbb{R}^{m\times n}$, \ahmadreza{where} $\uncertset$ is a convex and compact \ahmadreza{set in the form \eqref{Eq:uncertainty set form for slater} containing $\zero{m\times n}$ as a Slater point}.
 Here, we list all assumptions on the constraints in the forms (\ref{convex uncertainty convex quadratic}) and (\ref{convex uncertainty conic quadratic}) that we will make in this section, and use some of them in each theorem.
 \\
 \textbf{Assumption}:
 \begin{enumerate}[(A)]
 	\item \label{assumption: upper bound on the uncertainty} there exists $\Omega>0$ such that $\Vert \Delta\Vert_{2,2}\leq \Omega$ for all $\Delta\in \uncertset$.
 	\item \label{assumption: positive semidefinitness of a part} $\bar{ A}^T\bar{A}+2\Delta ^T\bar{ A}$ is positive semi-definite for all $\Delta \in \uncertset$.
 \end{enumerate}
 
 \ahmadreza{In addition to the discussion in Appendix \ref{APP: some sets for assumption psd for all} on how we can check these assumptions for a general uncertainty set $\uncertset$, here we show how to check them on two typically used uncertainty sets.
 	\\ \textbf{Ellipsoidal uncertainty set}: Let us assume that  $\mathcal{Z}=\left\{\Delta\in\mathbb{R}^{m\times n}:\hs \Vert\Delta\Vert_{2,2}\leq \rho \right\},$ for some $\rho>0$. Clearly, Assumption \eqref{assumption: upper bound on the uncertainty}  holds with $\Omega=\rho$. Furthermore, Assumption \eqref{assumption: positive semidefinitness of a part} can be checked using the following proposition:
 	\ahmadreza{\begin{proposition}\label{Prop:assumption psd for all, ellipsoid}
 			Let us assume that $\bar{ A}\in\mathbb{R}^{m\times n}$ is full column-rank and $\uncertset=\left\{\Delta\in\mathbb{R}^{m\times n}:\hs \Vert \Delta\Vert_{2,2}\leq \rho \right\}$. Then, Assumption \eqref{assumption: positive semidefinitness of a part} holds if and only if $\lambda_{\min}(\bar{ A}^T\bar{ A})\geq 4\rho^2,$ where $\lambda_{\min}(\bar{ A}^T\bar{ A})$ denotes the smallest eigenvalue of $\bar{ A}^T\bar{ A}.$
 		\end{proposition}
 		\begin{proof}{Proof.}
 			\ahmadreza{ 			Using Theorem 8.2.3 in \cite{bental2009robust},  Assumption \eqref{assumption: positive semidefinitness of a part} holds if and only if there exists a positive $\lambda$ such that
 				$$
 				\left[
 				\begin{aligned}
 				\lambda \identity{n}&\hspace{1cm}\rho\bar{ A}\\
 				\rho \bar{ A}^T& \hspace{0.5cm} \bar{ A}^T\bar{ A}-\lambda\identity{n}
 				\end{aligned}
 				\right]\succeq \zero{2n\times 2n}.
 				$$
 				Using the Schur Complement lemma (see, e.g., Appendix A.5.5 in  \cite{boyd2004convex}) the above linear matrix inequality is equivalent to $(\lambda-\rho^2)\bar{ A}^T\bar{ A}\succeq\lambda^2\identity{n}$. Hence, Assumption \eqref{assumption: positive semidefinitness of a part} holds if and only if 
 				\begin{eqnarray}
 				&\exists \lambda>\rho^2:\hs \bar{ A}^T\bar{ A}\succeq\frac{\lambda^2}{\lambda-\rho^2}\identity{n}\nonumber \\
 				\Leftrightarrow	&\hs	\exists \lambda>\rho^2:\hs \lambda_{\min}(\bar{ A}^T\bar{ A})\geq\frac{\lambda^2}{\lambda-\rho^2}\nonumber\\
 				\Leftrightarrow	&\hs\lambda_{\min}(\bar{ A}^T\bar{ A}) \geq 4\rho^2,\nonumber
 				\end{eqnarray}
 				where $\lambda_{\min}(\bar{ A}^T\bar{ A})$ denotes the smallest eigenvalue of $\bar{ A}^T\bar{ A}$, and where the last equivalence holds since $\frac{\lambda^2}{\lambda-\rho^2}$ is convex in $\lambda$ (for $\lambda>\rho^2$) with the minimum value of $4\rho^2.$}
 		\end{proof}
 	}
 	\noindent\textbf{Box uncertainty set}: Let us assume that  $\mathcal{Z}=\left\{\Delta\in\mathbb{R}^{m\times n}:\hs \Vert\Delta\Vert_\infty\leq \rho \right\},$ for some $\rho>0$. Using the following proposition, one can see that Assumption \eqref{assumption: upper bound on the uncertainty} holds with $\Omega=\rho\sqrt{nm}.$ 
 	\begin{proposition}\label{proposition: Omega for box uncertainty}
 		Let  $
 		\uncertset=\left\{\Delta\in \mathbb{R}^{\ahmadreza{m}\times n}: \;\Vert\Delta\Vert_\infty\leq \rho \right\}
 		$. Then,	$\sup_{\Delta\in \uncertset}\Vert\Delta\Vert_{2,2}=\rho\sqrt{nm}$.
 	\end{proposition} 
 	\begin{proof}{Proof.}
 		This follows directly from the definition of $\ell_\infty$ and spectral norms.
 	\end{proof}
 	Moreover, using Proposition \ref{Remark: approximation of assumption C} in Appendix \ref{APP: some sets for assumption psd for all}, Assumption \eqref{assumption: positive semidefinitness of a part} holds if 
 	\begin{equation}\label{Eq:approximation of DC formulation of assumption psd for all}
 	\min_{y\in\mathbb{R}^n}\left\{y^T\bar{ A}^T\bar{ A}y-2\rho \mathbb{1}^T \bar{ A} y:\hs\bar{ A} y\geq0,\hs \Vert y\Vert_1\leq1\right\} \geq0,
 	\end{equation}
 	where $\mathbb{1}\in\mathbb{R}^m$ is a vector whose components are all one. 
 } 
 \begin{remark}\label{Remark:simplicial uncertainty set is easy}
 	\ahmadreza{	Notice that if the uncertainty set is $\uncertset=\left\{\Delta\in\mathbb{R}^{m\times n}:\hs \Vert\Delta\Vert_1\leq\rho \right\}$, for some $\rho>0$, then \eqref{convex uncertainty} can be reformulated to a system of $2mn$ deterministic (conic) quadratic constraints because the uncertainty set contains $2mn$ vertices.} 
 \end{remark}
 Now, we proceed to the main results of this section. The following theorem provides tractable inner approximations of the constraints in the forms (\ref{convex uncertainty}) by replacing the quadratic term in the uncertain parameter with a linear upper bound.
 
 \begin{theorem}\label{Th: inner approximation} Let \ahmadreza{Assumption (\ref{assumption: upper bound on the uncertainty}) hold}. Then:
 	\begin{enumerate}[(I)]
 		\item \label{Th: inner approximation, convex quadratic} $y\in \mathbb{R}^n$ satisfies (\ref{convex uncertainty convex quadratic}) if there exists $\uncertA \in \mathbb{R}^{n \times n}$ satisfying the convex system
 		\begin{equation}\label{final robust counterpart of perturbation by vital}
 		\begin{matrix}
 		\trace{( \bar{ A}^T\bar{A}+\Omega^2\identity{n})}{\uncertA}+\delta^*_{\uncertset}(2\bar{ A}\uncertA +D^T\adj a^T)+\bar{b}^Ty+c\leq 0,&\hs
 		\left[\begin{matrix}
 		\uncertA &\adj\\
 		\adj^T&1
 		\end{matrix}\right]\succeq \zero{n+1\times n+1}.
 		\end{matrix}
 		\end{equation}
 		\item \label{Th: inner approximation, conic quadratic} $y\in \mathbb{R}^n$ satisfies (\ref{convex uncertainty conic quadratic}) if there exist $\uncertA \in \mathbb{R}^{n \times n}$ and $\scalervar \in \mathbb{R}$ satisfying the convex system
 		\begin{equation}\label{final robust counterpart of perturbation conic by vital}
 		\begin{matrix}
 		\trace{( \bar{ A}^T\bar{A}+\Omega^2\identity{n})}{\uncertA}+\delta^*_{\uncertset}(2\bar{A}\uncertA+D^T\adj a^T)+\bar{ b}^T\adj +c+\frac{\scalervar}{4}\leq 0,\hs
 		\left[\begin{matrix} \uncertA&\adj\\\adj^T&\scalervar
 		\end{matrix}\right]\succeq \zero{n+1\times n+1}.
 		\end{matrix}
 		\end{equation}
 	\end{enumerate}
 \end{theorem}
 \begin{proof}{Proof.}
 	(\ref{Th: inner approximation, convex quadratic})  $y\in \mathbb{R}^n$ satisfies (\ref{convex uncertainty convex quadratic}) if and only if 
 	\begin{equation}\label{robust counterpart of perturbation}
 	y^T\bar{ A}^T\bar{ A}y+2y^T\bar{A}^T\Delta  y+\Eucledianorm{ \Delta \adj}^2+(D\Delta a)^Ty+\bar{ b}^Ty+c\leq 0, \hspace{.3cm }\forall \Delta \in \uncertset. 
 	\end{equation}
 	Replacing $\Eucledianorm{ \Delta \adj}^2$ by its upper bound $\Omega^2\Eucledianorm{\adj}^2$ implies that $y\in \mathbb{R}^n$ satisfies (\ref{robust counterpart of perturbation}) if it satisfies
 	\begin{equation}\label{inner approximation of the robust counterpart of perturbation}
 	y^T\bar{ A}^T\bar{ A}y+2y^T\bar{ A}^T\Delta  y+\Omega^2y^T y+(D\Delta a)^Ty+\bar{ b}^Ty+c\leq 0, \hspace{.3cm }\forall \Delta \in \uncertset.
 	\end{equation}
 	Setting $\mathcal{U}=\{(\bar{ A}^T\bar{A}+2\bar{ A}^T\Delta  +\Omega^2 \identity{n}, D\Delta a): \hspace{.2cm}\Delta \in \uncertset\}$, (\ref{inner approximation of the robust counterpart of perturbation}) is equivalent to
 	$$
 	y^TBy+(\bar{b}+d)^Ty+c\leq 0 \hspace{.3cm }\forall (B,d)\in \mathcal{U}.
 	$$
 	For any $(B,d)\in \mathcal{U}$, $B$ is positive semi-definite since 
 	$B=(\bar{A}+\Delta)^T(\bar{A}+\Delta)+\Omega^2 \identity{n}-\Delta^T\Delta\succeq 0_{n\times n}.$
 	%
 	So, by applying Theorem \ref{Th: RC for conic and convex quadratic}(\ref{Th: item, convex quadratic}) and Lemma \ref{Lemma: composition support}(\ref{general format of affinity}), $y\in \mathbb{R}^n$ satisfies (\ref{robust counterpart of perturbation}) if there exists $\uncertA \in \mathbb{R}^{n\times n}$ such that $\adj$ and $\uncertA$ satisfy (\ref{final robust counterpart of perturbation by vital}).
 	\\
 	(\ref{Th: inner approximation, conic quadratic}) The proof is similar to part (I).
 \end{proof}
 
 
 In the next theorem we derive tractable outer approximations of the constraints in the forms (\ref{convex uncertainty}).
 
 \begin{theorem}\label{Th: tractable outer approxmation}
 	Let Assumption (\ref{assumption: positive semidefinitness of a part}) holds. Then:
 	\begin{enumerate}[(I)]
 		\item \label{Th: tractable outer approximation item convex quadratic } if	$y$ satisfies (\ref{convex uncertainty convex quadratic}), then there exists $\uncertA\in \mathbb{R}^{n\times n}$ satisfying the convex system
 		\begin{equation}\label{Eq: tractable outer approximation convex quadratic}  
 		\begin{matrix}
 		\trace{\bar{ A}^T\bar{ A}}{\uncertA}+\delta^*_{\uncertset}(2\bar{ A}\uncertA+D^T\adj a^T)+\bar{ b}^Ty+c\leq 0,&\hs\hs\hs
 		\left[\begin{matrix}
 		\uncertA &\adj\\
 		\adj^T&1
 		\end{matrix}\right]\succeq \zero{n+1\times n+1}. 
 		\end{matrix}
 		\end{equation}
 		\item \label{Th: outer approximation, conic quadratic} if $y\in \mathbb{R}^n$ satisfies (\ref{convex uncertainty conic quadratic}), then there exist $\uncertA \in \mathbb{R}^{n \times n}$ and $\scalervar \in \mathbb{R}$ satisfying the convex system
 		\begin{equation}\label{Eq: tractable outer approximation conic quadratic}
 		\begin{matrix}
 		\trace{\bar{ A}^T\bar{ A}}{\uncertA}+\delta^*_{\uncertset}(2\bar{ A}\uncertA +D^T\adj a^T)+\bar{ b}^T\adj +c+\frac{\scalervar}{4}\leq 0,&\hs
 		\left[\begin{matrix} \uncertA&\adj\\\adj^T&\scalervar
 		\end{matrix}\right]\succeq \zero{n+1\times n+1}.
 		\end{matrix}
 		\end{equation} 
 	\end{enumerate}
 \end{theorem}
 \begin{proof}{Proof.}
 	(\ref{Th: tractable outer approximation item convex quadratic }) It is clear that $\adj$ satisfies (\ref{convex uncertainty convex quadratic}) if and only if $\adj$ satisfies (\ref{robust counterpart of perturbation}). Replacing $\Eucledianorm{\Delta \adj}^2$ with its lower bound $0$ implies that if $y\in \mathbb{R}^n$ satisfies (\ref{robust counterpart of perturbation}) then 
 	\begin{equation}\label{Eq: deleting the quadratic uncertainty}
 	y^T\bar{ A}^T\bar{ A}y+2y^T\bar{ A}^T\Delta  y+(D\Delta a)^Ty+\bar{ b}^Ty+c\leq 0, \hspace{.3cm }\forall \Delta \in \uncertset.
 	\end{equation}
 	Setting $\mathcal{U}=\{(\bar{ A}^T\bar{A}+2\bar{ A}^T\Delta  , D\Delta a): \hspace{.2cm}\Delta \in \uncertset\}$, and using Theorem \ref{Th: RC for conic and convex quadratic}(\ref{Th: item, convex quadratic}) and Lemma \ref{Lemma: composition support}(\ref{general format of affinity}) completes the proof.
 	\\
 	(\ref{Th: outer approximation, conic quadratic}) The proof is similar to the previous part. 
 \end{proof}
 
 In the next theorem we provide an upper bound on the violation errors of (\ref{convex uncertainty convex quadratic}) and (\ref{convex uncertainty conic quadratic}) for the solutions that satisfy the outer approximations
 (\ref{Eq: tractable outer approximation convex quadratic}) and (\ref{Eq: tractable outer approximation conic quadratic}), respectively.
 
 \begin{theorem}\label{Th: upper bound on the violation error of the outer approximation}
 	Let Assumptions (\ref{assumption: upper bound on the uncertainty}) and (\ref{assumption: positive semidefinitness of a part}) hold. Then, 
 	\begin{enumerate}[(I)]
 		\item \label{Th: upper bound on the violation error of the outer approximation convex quadratic} if $\adj\in \mathbb{R}^n$ and $\uncertA\in \mathbb{R}^{n \times n}$ satisfy (\ref{Eq: tractable outer approximation convex quadratic}), then $\adj$ violates (\ref{convex uncertainty convex quadratic}) by at most $\Omega^2\Eucledianorm{ \adj}^2$.
 		\item \label{Th: upper bound on the violation error of the outer approximation conic quadratic} if $\adj\in \mathbb{R}^n$ and $\uncertA\in \mathbb{R}^{n \times n}$ satisfy (\ref{Eq: tractable outer approximation conic quadratic}), then $\adj$ violates (\ref{convex uncertainty conic quadratic}) by at most $\Omega\Eucledianorm{\adj}$.
 	\end{enumerate}
 \end{theorem}
 \begin{proof}{Proof.}
 	(\ref{Th: upper bound on the violation error of the outer approximation convex quadratic})  Let $\adj\in \mathbb{R}^n$ and $\uncertA\in \mathbb{R}^{n \times n}$ satisfy (\ref{Eq: tractable outer approximation convex quadratic}). Then, $\adj$ satisfies (\ref{Eq: deleting the quadratic uncertainty}). Therefore, 
 	\begin{equation}\label{Eq: deleting the quadratic term in uncertainty with max}
 	\max _{\Delta\in \uncertset}\{y^T\bar{ A}^T\bar{ A}y+2y^T\bar{ A}^T\Delta  y+(D\Delta a)^Ty+\bar{b}^Ty+c\}\leq 0.
 	\end{equation}
 	As it is mentioned in the proof of Theorem \ref{Th: inner approximation}(\ref{Th: inner approximation, convex quadratic}), (\ref{convex uncertainty convex quadratic}) is equivalent to (\ref{robust counterpart of perturbation}). Therefore, we have
 	$$
 	\begin{aligned}
 	&\;\max _{\Delta\in \uncertset}\{y^T\bar{ A}^T\bar{ A}y+2y^T\bar{ A}^T\Delta  y+\Eucledianorm{\Delta \adj}^2+(D\Delta a)^Ty+\bar{ b}^Ty+c\}\\
 	&\;\leq y^T\bar{ A}^T\bar{ A}y+\bar{ b}^Ty+c+\max _{\Delta\in \uncertset}\{2y^T\bar{ A}^T\Delta  y+(D\Delta a)^Ty\}+\max _{\Delta\in \uncertset}\Eucledianorm{\Delta \adj}^2\\
 	&\;\leq y^T\bar{ A}^T\bar{ A}y+\bar{ b}^Ty+c+\max _{\Delta\in \uncertset}\{2y^T\bar{ A}^T\Delta  y+(D\Delta a)^Ty\}+\Omega^2\Eucledianorm{ \adj}^2
 	\; \leq \Omega^2\Eucledianorm{  \adj}^2,
 	\end{aligned}
 	$$ 
 	where the last inequality follows from (\ref{Eq: deleting the quadratic term in uncertainty with max}).
 	\\
 	(\ref{Th: upper bound on the violation error of the outer approximation conic quadratic}) It is clear that (\ref{convex uncertainty conic quadratic}) is equivalent to 
 	$$
 	\sqrt{y^T\bar{ A}^T\bar{ A}y+2y^T\bar{ A}^T\Delta  y+\Eucledianorm{ \Delta \adj}^2+(D\Delta a)^Ty}+\bar{ b}^Ty+c\leq 0, \hspace{.3cm }\forall \Delta \in \uncertset.
 	$$
 	Similar to the previous part, if $\adj$ comes from the outer approximation (\ref{Eq: tractable outer approximation conic quadratic}), then we have
 	$$
 	\begin{aligned}
 	&\;\sqrt{\max _{\Delta\in \uncertset}\{y^T\bar{ A}^T\bar{ A}y+2y^T\bar{ A}^T\Delta  y+\Eucledianorm{\Delta \adj}^2+(D\Delta a)^Ty\}}+\bar{ b}^Ty+c\\
 	&\;\leq \sqrt{y^T\bar{ A}^T\bar{ A}y+\max _{\Delta\in \uncertset}\{2y^T\bar{ A}^T\Delta  y+(D\Delta a)^Ty\}}+\max _{\Delta\in \uncertset}\Eucledianorm{\Delta \adj}+\bar{ b}^Ty+c\\
 	&\;\leq \sqrt{y^T\bar{ A}^T\bar{ A}y+\max _{\Delta\in \uncertset}\{2y^T\bar{ A}^T\Delta  y+(D\Delta a)^Ty\}}+\Omega\Eucledianorm{ \adj}+\bar{ b}^Ty+c
 	\; \leq \Omega\Eucledianorm{  \adj},
 	\end{aligned}
 	$$ 
 	where the first inequality holds because of the fact that $\sqrt{f+g}\leq \sqrt{f}+\sqrt{g}$ for any $f,g\geq 0$.
 \end{proof}
 
 \begin{remark}
 	\ahmadreza{Until} now, we have considered problems containing uncertainties in their constraint parameters. This is without loss of generality, since if we have a problem with uncertainty in the parameters of the objective function, then we can use the epigraph formulation to shift the uncertainty to a constraint. \hfill \qed
 \end{remark}


 \section{Data-driven uncertainty set}\label{Sec: new statistical uncertainty set}
 A usual way of constructing an uncertainty set is by using historical data and statistical tools, such as hypothesis testing (\cite{Bertsimas2017data}), or asymptotic confidence sets (\cite{ben2013robust}).  In this section, we  use the latter to design an uncertainty set for a vector consisting of the mean and vectorized covariance matrix.

 For notational simplicity, we explain how to construct an uncertainty set for the two dimensional case; the extension to higher dimensions is straightforward.   For the two dimensional case, assume that $\left(x \atop z\right)$ is a random vector with components $x, z$ and set $\mu_x=\mathbb{E}(x)$, $\mu_z=\mathbb{E}(z)$, $\sigma^2_x=\mathbb{E}(x-\mu_x)^2$, $\sigma^2_z=\mathbb{E}(z-\mu_z)^2$, $\sigma_{xz}=\mathbb{E}(x-\mu_x)(z-\mu_z)$, and $\mu_{kl}=\mathbb{E}(x-\mu_x)^k(z-\mu_z)^l$, $k,l=0,1,2,...$. Assume that the fourth moments exist, which means that $\mu_{kl}$ exists when $k+l\leq 4$, $k,l=0,1,2,3,4$. This assumption can be tested using the result in \cite{trapani2016testing}. Now, consider a random sample of size $n$, $\left(x_i \atop z_i\right)$, $i=1,...,n$. Set 
 $$
 \bar{x}=\frac{1}{n}\sum_{i=1}^nx_i,\hspace{0.2cm }\bar{z}=\frac{1}{n}\sum_{i=1}^nz_i,\hspace{0.2cm }S^2_x=\frac{1}{n}\sum_{i=1}^n(x_i-\bar{x})^2,
 \hspace{0.2cm }S^2_z=\frac{1}{n}\sum_{i=1}^n(z_i-\bar{z})^2,\hspace{0.2cm }S_{xz}=\frac{1}{n}\sum_{i=1}^n(x_i-\bar{x})(z_i-\bar{z}).
 $$
 Using the Central Limit Theorem (Example 2.18 in \cite{van2000asymptotic}) and the Delta Method (Theorem 3.1 in \cite{van2000asymptotic}), and setting 
 $$
 \mathcal{Y}=
 \left(
 \mu_x, \ 
 \mu_z, \ 
 \mathbb{E}(x^2),\ 
 \mathbb{E}(xz),\ 
 \mathbb{E}(z^2)
 \right)^T
 ,\hspace{0.3cm}
 Y_n=
 \left(
 \bar{x},  \ 
 \bar{z}, \ 
 \frac{1}{n}\sum_{i=1}^nx_i^2,\ 
 \frac{1}{n}\sum_{i=1}^nx_iz_i,\ 
 \frac{1}{n}\sum_{i=1}^nz_i^2
 \right)^T,$$
 it follows for any differentiable function $\phi:\mathbb{R}^5\rightarrow \mathbb{R}^m$ that $\sqrt{n}\left(\phi(Y_n)-\phi(\mathcal{Y})\right)$ converges in distribution to the normal distribution $N(0,\nabla\phi(\theta)\Sigma\nabla\phi(\theta)^T)$, where $\Sigma$ and $\nabla\phi$ are the covariance matrix of $\left(x,z,x^2,xz,z^2\right)^T-\mathcal{Y}$ and the Jacobian matrix of $\phi$, respectively. Letting $$\phi(x_1,...,x_5)=\left(
 x_1,\ x_2, \ 
 x_3-x_1^2,\ 
 x_4-x_1x_2,\ 
 x_5-x_2^2
 \right)^T,$$
 it is easy to show, similar to Example 3.2 in \cite{van2000asymptotic}, that
 \begin{equation}\label{Eq: distribution of mu and sigma}
 \sqrt{n}\left(T_n-\theta\right)\xrightarrow[n\rightarrow \infty]{d} N\left(
 \left(
 \begin{matrix}
 0 \\
 0 \\
 0\\
 0\\
 0\\
 \end{matrix}
 \right),
 \underbrace{\left(
 	\begin{matrix}
 	\mu_{20}&\mu_{11}&\mu_{30}&\mu_{21}&\mu_{12}\\
 	\mu_{11}&\mu_{02}&\mu_{21}&\mu_{12}&\mu_{03}\\
 	\mu_{30}&\mu_{21}&\mu_{40}-\mu_{20}^2&\mu_{31}-\mu_{11}\mu_{20}&\mu_{22}-\mu_{20}\mu_{02}\\
 	\mu_{21}&\mu_{12}&\mu_{31}-\mu_{11}\mu_{20}&\mu_{22}-\mu_{11}^2&\mu_{13}-\mu_{11}\mu_{02}\\
 	\mu_{12}&\mu_{03}&\mu_{22}-\mu_{20}\mu_{02}&\mu_{13}-\mu_{11}\mu_{02}&\mu_{04}-\mu_{02}^2
 	\end{matrix}
 	\right)}_{V}
 \right),
 \end{equation}
 where
 $
 \theta=\phi(\mathcal{Y})=
 \left(
 \mu_x, \ 
 \mu_z, \ 
 \sigma^2_x,\ 
 \sigma_{xz},\ 
 \sigma^2_z
 \right)^T
 $, $
 T_n=\phi(\mathcal{Y}_n)=
 \left(
 \bar{x}, \ 
 \bar{z}, \ 
 S^2_x,\ 
 S_{xz},\ 
 S^2_z
 \right)^T,$
 and  $\xrightarrow[n\rightarrow \infty]{d}$ means convergence in distribution when the size of the random sample goes to infinity.
 
 Let $\hat{V}$ and $\hat{\theta}$ be consistent estimates of $V$ and $\theta$ defined in (\ref{Eq: distribution of mu and sigma}), respectively.
 Then, asymptotically  with $(1-\alpha)\%$ confidence, $\theta$ belongs to the following ellipsoid:
 $$
 \mathcal{U}:=\left\{\theta:\hspace{0.2cm} n\left(\hat{\theta}-\theta\right)^T\hat{V}^{-1}\left(\hat{\theta}-\theta\right)\leq \chi^2_{rank(V),1-\alpha} \right\},
 $$ 
 where $\chi^2_{d, 1-\alpha}$ denotes the $(1-\alpha)$ percentile of the Chi-square distribution with $d$ degrees of freedom.
 
 To use the results of Section \ref{Sec:exact formulation}, we reformulate the uncertainty set $\mathcal{U}$. 
 Setting 
 \begin{equation}\label{Eq:notation for A and mu and sigma}
 \begin{aligned}
 &\Psi=\left[\Psi_\mu\atop \Psi_\Sigma\right],\; \;\Psi_\mu=\left(\begin{matrix}
 1&0&0&0&0\\
 0&1&0&0&0\\
 \end{matrix}\right),\;\;\Psi_\Sigma=\left(\begin{matrix}
 0&0&1&0&0\\
 0&0&0&\sqrt{2}&0\\
 0&0&0&0&1\\
 \end{matrix}\right),
 &\mu=[\mu_x  \ \mu_z]^T,\;\; \Sigma=\left[\begin{matrix}
 \sigma^2_x&\sigma_{xz}\\\sigma_{xz}& \sigma^2_z,
 \end{matrix}\right],
 \end{aligned}\end{equation}
 due to positive semi-definiteness of $\Sigma$, with $(1-\alpha)\%$ confidence
 $$
 \left(
 \begin{matrix}
 \mu\\svec(\Sigma)
 \end{matrix}
 \right)\in \bar{\mathcal{U}}:=\hat{\mathcal{U}}\cap
 \left\{\gamma:\hspace{0.2cm} n\left(\Psi\hat{\theta}-\gamma\right)^T\Psi^{-1}\hat{V}^{-1}\Psi^{-1}\left(\Psi\hat{\theta}-\gamma\right)\leq \chi^2_{rank(V),1-\alpha} \right\},
 $$ 
 where $\hat{\mathcal{U}}=\left\{\left(\gamma_\mu\atop \gamma_\Sigma \right):\;\gamma_\Sigma=svec(M),\;M\succeq\zero{n\times n} \right\}$. 
 Letting $R^TR$ be the Cholesky factorization of $\hat{V}^{-1}$, i.e., $\hat{V}^{-1}=R^TR$, $\bar{\mathcal{U}}$ can be rewritten as 
 $$
 \begin{aligned}
 \bar{\mathcal{U}}=&\hs\hat{\mathcal{U}}\cap\left\{\gamma:\hs\Eucledianorm{R\Psi^{-1}\left(\gamma -\Psi\hat{\theta}\right)}\leq \sqrt{\frac{\chi^2_{rank(V),1-\alpha}}{n}}  \right\}\\
 =&\hs\hat{\mathcal{U}}\cap\left\{\Psi R^{-1}\nu+\Psi\hat{\theta}:\hs\Vert\nu\Vert _2\leq \sqrt{\frac{\chi^2_{rank(V),1-\alpha}}{n}}  \right\}.
 \end{aligned}
 $$
 Hence, by letting the estimated mean vector and covariance matrix based on the random sample be $\hat{\mu}$ and $\hat{\Sigma}$, respectively, we have
 \begin{equation}\label{Eq: the new statistical uncertainty set}
 \bar{\mathcal{U}}=\hat{\mathcal{U}}\cap\left\{\Psi R^{-1}\nu+\left(\hat{\mu} \atop svec(\hat{\Sigma}) \right):\hs \Vert\nu\Vert _2\leq \sqrt{\frac{\chi^2_{rank(V),1-\alpha}}{n}}  \right\}.
 \end{equation}

 \begin{remark}\label{Remark: invertibility of V}
 	If $V$ is not invertible, then one can use a generalized inverse, such as the \textit{Moore-Penrose} inverse.	
 	\hfill \qed
 \end{remark}
 \begin{remark}\label{Remark: generalization to higher dimension}
 	The construction of the  uncertainty set can straightforwardly be  extended to higher dimensions  using suitable $\phi,\ \Psi,$ and $V$. Details are omitted for brevity of exposition. 
 	\hfill \qed
 \end{remark}
 \begin{remark}\label{Remark:time-series}
 	\ahmadreza{The uncertainty set $\bar{ \mathcal{U}}$ is constructed for a random sample. Analogously, one can construct an uncertainty set for a time-series under appropriate assumptions; see, e.g., Section 2.2 in the book \cite{financialEco}. }\hfill\qed
 \end{remark}
 Now, consider a convex quadratic constraint 
 \begin{equation}\label{Eq: a general quadratic const. with mu and sigma}
 \adj^T\Sigma\adj+\mu^T\adj+c\leq 0,
 \end{equation} 
 where $\mu$ and $\Sigma$ are the mean vector and covariance matrix of a random vector. By using the uncertainty set $\bar{\mathcal{U}}$ in (\ref{Eq: the new statistical uncertainty set}) and Example \ref{Ex:general statical uncertainty set}, the RC of (\ref{Eq: a general quadratic const. with mu and sigma}) is
 \begin{equation}\label{Eq: robust counterpart of a quadratic constraint with sigma and mu}
 \begin{aligned}
 &\hat{\mu}^T\adj+\trace{\hat{\Sigma}}{\uncertA}+\rho\Eucledianorm{\left(\Psi R^{-1}\right)^T\left(\adj \atop svec(\uncertA)\right)}+c\leq 0,\;
 &\left[\begin{matrix}
 \uncertA & \adj\\
 \adj^T&1
 \end{matrix}\right]\succeq \zero{n+1\times n+1},
 \end{aligned}
 \end{equation}
 where $\rho=\sqrt{\frac{\chi^2_{rank(V),1-\alpha}}{n}}$.
 
 Let $\mu$ and $\Sigma$ be the actual population mean vector and covariance matrix, respectively. Then, $\theta=\left(\mu\atop svec(\Sigma)\right)$ belongs to the uncertainty set $\bar{ \mathcal{U}}$ asymptotically  with confidence level $(1-\alpha)\%$. This, roughly speaking, means that the uncertainty set not only contains $\theta$ but also many more points. 
 Therefore, $\adj$ that satisfies (\ref{Eq: robust counterpart of a quadratic constraint with sigma and mu}) is \ahmadreza{asymptotically} immunized against some extra  $\mu$ and $\Sigma$ and hence conservative. 
 
 Another way of dealing with the uncertainty in $\theta$ is by making use of the chance constraint 
 $
 Prob\left(\adj^T\Sigma\adj+\mu^T\adj+c\leq 0\right)\geq 1-\alpha,
 $
 where $\alpha>0$ is close to $0$.  In what follows, we elaborate more on this chance constraint and provide a reformulation and relaxation of it.
 
 For any vector $\beta$, (\ref{Eq: distribution of mu and sigma}) implies that
 $
 \sqrt{n}\left(\beta ^TT_n-\beta^T\theta\right)\xrightarrow[n\rightarrow \infty]{d} N\left(0,\beta^TV\beta\right).
 $
 By setting $\beta=\Psi\left(\adj \atop  svec(\adj\adj^T)\right)$, it follows straightforwardly that the (asymptotic) chance constraint with probability of $1-\alpha$ is equivalent to 
 \begin{equation}\label{Eq: reformulated chance constraint}
 \frac{z_{1-\alpha}}{\sqrt{n}}\sqrt{\beta^T\hat{V}\beta}+\hat{\mu}^T\adj+\adj^T\hat{\Sigma}\adj+c\leq 0,
 \end{equation}
 where $z_{1-\alpha}$ is the $1-\alpha$ percentile of the standard normal distribution. 
 Clearly (\ref{Eq: reformulated chance constraint}) is equivalent to the set of constraints
 $$
 \frac{z_{1-\alpha}}{\sqrt{n}}\Vert R^{T^{-1}}\beta\Vert+\hat{\mu}^T\adj+\adj^T\hat{\Sigma}\adj+c\leq 0, \;\;\beta=\Psi\left(\adj \atop  svec(W)\right),\;\;W=\adj\adj^T,
 $$
 where $R$ is the Cholesky factorization of $\hat{V}^{-1}$. The constraint $W=\adj\adj^T$ is nonconvex, so we relax it to $W\succeq\adj\adj^T$, which is a semi-definite representable constraint. Hence, 
 \begin{equation}\label{Eq: relaxed chance constraint reformulation}
 \begin{aligned}
 \frac{z_{1-\alpha}}{\sqrt{n}}\Vert R^{T^{-1}}\beta\Vert+\hat{\mu}^T\adj+\adj^T\hat{\Sigma}\adj+c\leq 0, &\;\;\beta=\Psi\left(\adj \atop  svec(W)\right),&
 \left[\begin{matrix}
 \uncertA &\adj\\
 \adj^T&1
 \end{matrix}\right]\succeq \zero{n+1\times n+1},
 \end{aligned}
 \end{equation}
 is a relaxation of (\ref{Eq: reformulated chance constraint}). In the next proposition, we provide a relation between solutions that satisfy \eqref{Eq: robust counterpart of a quadratic constraint with sigma and mu} and the ones satisfying \eqref{Eq: relaxed chance constraint reformulation}. 
 
 \begin{proposition}\label{Proposition: chance relaxation is better than robust}
 	Let $\left(\adj, W\right)$ be a solution that satisfies \eqref{Eq: robust counterpart of a quadratic constraint with sigma and mu}. Then $\left(\adj, W\right)$ also satisfies \eqref{Eq: relaxed chance constraint reformulation}.
 \end{proposition}
 \begin{proof}{Proof.}
 	Appendix \ref{Proof: proposition chance relaxation is better than robust}.
 \end{proof}
 
 Even though Proposition \ref{Proposition: chance relaxation is better than robust} asserts that \eqref{Eq: robust counterpart of a quadratic constraint with sigma and mu} is more conservative than \eqref{Eq: relaxed chance constraint reformulation}, we cannot conclude that \eqref{Eq: robust counterpart of a quadratic constraint with sigma and mu} is more conservative than \eqref{Eq: reformulated chance constraint}. This is because $\uncertA=\adj\adj^T$ is not necessarily satisfied for solutions of \eqref{Eq: relaxed chance constraint reformulation}.
 
 \begin{remark}
 	After solving the problem containing (\ref{Eq: relaxed chance constraint reformulation}),  
 	if $W=\adj\adj^T$ is not satisfied, then $\adj$ is not feasible for \eqref{Eq: reformulated chance constraint}. However, $\adj$ is strictly feasible for 
 	\begin{equation}\label{Eq:portfolio onjective}
 	\adj^T\Sigma \adj+\mu^T\adj+c\leq  0
 	\end{equation}
 	for the nominal scenario $(\hat{\mu},\hat{\Sigma})$. Therefore, $\adj$ may also be feasible for other scenarios and hence is more robust than the nominal solution. 
 	\hfill \qed \end{remark}
 
 %
 Consider a solution $\bar{\adj}$ that satisfies (\ref{Eq: robust counterpart of a quadratic constraint with sigma and mu}) where the uncertainty set is constructed using the desired confidence level $1-\bar{ \alpha}$. For this solution,	according to the above discussion, $Prob\left(\bar{\adj}^T\Sigma\bar{\adj}+\mu^T\bar{\adj}+c\leq 0\right)$ might be larger than the desired confidence level $1-\bar{ \alpha}$. If so, then by decreasing the confidence level that is used in the construction of the uncertainty set and considering $\tilde{\adj}$ that satisfies (\ref{Eq: robust counterpart of a quadratic constraint with sigma and mu}), $Prob\left(\tilde{\adj}^T\Sigma\tilde{\adj}+\mu^T\tilde{\adj}+c\leq 0\right)$ gets closer to the desired confidence level $1-\bar{ \alpha}$. In our numerical experiments, we check for different instances which confidence level should be used in the construction of the uncertainty set such that for the robust solution the constraint $\adj^T\Sigma\adj+\mu^T\adj+c\leq0$ is satisfied with probability close to the desired confidence level.
 \section{Applications}\label{Sec: numerical result}
 In this section, we apply the results of the previous sections to a robust portfolio choice, norm approximation, and regression line problem. \ahmadreza{All computations in this paper were carried out with MATLAB 2016a using YALMIP \cite{Lofberg2004} to pass the optimization problems to MOSEK 8.1.0.80 \cite{mosek2015mosek}.  }
 \subsection{Mean-Variance portfolio problem}\label{Sec: mean-variance}
 In this subsection, we describe a formulation for a mean-variance portfolio problem (Chapter 2 in \cite{fabozzi2007robust}), and use the results of Section \ref{Sec: new statistical uncertainty set} to construct an uncertainty set and to derive a tractable reformulation of the robust counterpart. \\
 \textbf{Problem formulation:} We consider a mean-variance portfolio problem with $n$ assets. Let $\mu$ and $\Sigma$ be the expectation and covariance matrix of the return vector $r=(r_1,...,r_n)$, respectively. 
 One formulation of a  mean-variance portfolio problem is to model the trade-off between the risk and mean return in the objective function using a risk-aversion coefficient $\lambda$: 
 \begin{equation}\label{Risk aversion portfolio problem}
 \begin{aligned}
 \max_{\omega}\left\{  \mu^T\omega- \lambda \omega^T \Sigma \omega \;:
 \hspace{0.3cm} \mathbb{1}^T\omega=1,\;\;\omega\geq0\right\},
 \end{aligned}
 \end{equation}
 where $\mathbb{1}=[1,1,...,1]^T$. 
 The risk aversion coefficient is determined by the decision maker. When it is small, it means that the mean return is more important than the corresponding risk and it leads to a more risky portfolio than when the risk-aversion coefficient is large. 
 
 In practice $\mu$ and $\Sigma$ are typically estimated from a set of historical data, which makes them sensitive to sampling inaccuracy. There are several ways of defining uncertainty sets for the expected return vector and asset return covariance matrix, e.g., see Chapter 12 in  \cite{fabozzi2007robust}. In this section, we use $\bar{\mathcal{U}}$ defined in (\ref{Eq: the new statistical uncertainty set}), i.e., the uncertainty set constructed for $\left(\mu\atop svec(\Sigma)\right)$.
 Using (\ref{Eq: robust counterpart of a quadratic constraint with sigma and mu}), the robust counterpart of 
 (\ref{Risk aversion portfolio problem}) with uncertainty set $\bar{\mathcal{U}}$ reads
 \begin{equation}\label{robust risk aversion portfolio}
 \begin{aligned}
 \max_{\omega, \uncertA}&\hspace{.2cm} \hat{\mu}^T\omega-\lambda tr(\hat{\Sigma} \uncertA)-\rho \Eucledianorm{\left(\Psi R^{-1}\right)^T\left(-\omega\atop \lambda svec(\uncertA)\right)}\\
 \mbox{s.t.}&\left[\begin{matrix}
 \uncertA& \omega\\
 \omega^T& 1
 \end{matrix}\right] \succeq \zero{n+1\times n+1},\hspace{.2cm} \mathbb{1}^T\omega=1,
 \;\;\omega\geq0,
 \end{aligned}
 \end{equation}
 where $\rho=\sqrt{\frac{\chi^2_{rank(V),1-\alpha}}{n}}$, $\hat{\mu}$, $\Sigma$, and $\hat{V}$ are consistent estimates of $\mu$, $\Sigma$, and $V$, with $V$ and $\Psi$ as in (\ref{Eq: distribution of mu and sigma}) and (\ref{Eq:notation for A and mu and sigma}), respectively, but formulated for the higher dimensional case, and $R$ is the Cholesky factorization of $\hat{V}^{-1}$.
 
 Furthermore, by setting
 $
 \beta=\Psi\left(\omega \atop -\lambda svec(\uncertA)\right),
 $
 and using the relaxed chance constraint (\ref{Eq: relaxed chance constraint reformulation}), the robust counterpart of problem  (\ref{Risk aversion portfolio problem}) with confidence $(1-\alpha)\%$ is approximated by
 \begin{equation}\label{Eq: eqivalent risk aversion portfolio with chance constraint}
 \begin{aligned}
 \max_{\omega} &\hspace{0.3cm} \hat{\mu}^T\omega-\lambda\omega^T\hat{\Sigma}\omega-\frac{z_{1-\alpha}}{\sqrt{n}}\Eucledianorm{ R^{T^{-1}} \beta} \\
 \mbox{s.t.}&\hspace{0.3cm}     	\left[\begin{matrix}
 \uncertA &\omega\\
 \omega^T&1
 \end{matrix}\right]\succeq \zero{n+1\times n+1},\;\; \mathbb{1}^T\omega=1,\;\;\omega\geq0.
 \end{aligned}
 \end{equation}
 \\ \textbf{Numerical evaluation:} To evaluate the above robust counterparts, we use the monthly average value weighted return of 5 and 30 industries from 1956 \ahmadreza{until} 2015, obtained from  ``Industry Portfolios" data on the website \url{http://mba.tuck.dartmouth.edu/pages/faculty/ken.french/data_library.html}. The data are  monthly returns, but to present the results, we report the annualized returns (obtained by multiplying the expected monthly return by $12$) and the annualized risk (multiplication of the standard deviation by $\sqrt{12}$). 
 Furthermore, we set the risk aversion coefficient $\lambda$ to $3$.

 We have solved the following three problems: (\ref{Risk aversion portfolio problem}) with nominal values for $\mu$ and $\Sigma$ estimated from the data, which we call \textit{Nominal problem}; (\ref{robust risk aversion portfolio}), which we call \textit{Robust problem}; and (\ref{Eq: eqivalent risk aversion portfolio with chance constraint}), which we call \textit{Chance problem}, due to the chance constraint.
 
 We first check the behavior of $Prob(\mu^T\omega^*- \lambda \omega^{*^T} \Sigma \omega^*\geq z^*)$ as a function of the confidence level used to construct the uncertainty set, where $\omega^*$ and $z^*$ are the robust solution and corresponding robust objective value, respectively. 
 \begin{figure}
 	\begin{subfigure}{0.5\linewidth}
 		\centering
 		\includegraphics[width=8cm, height=5cm]{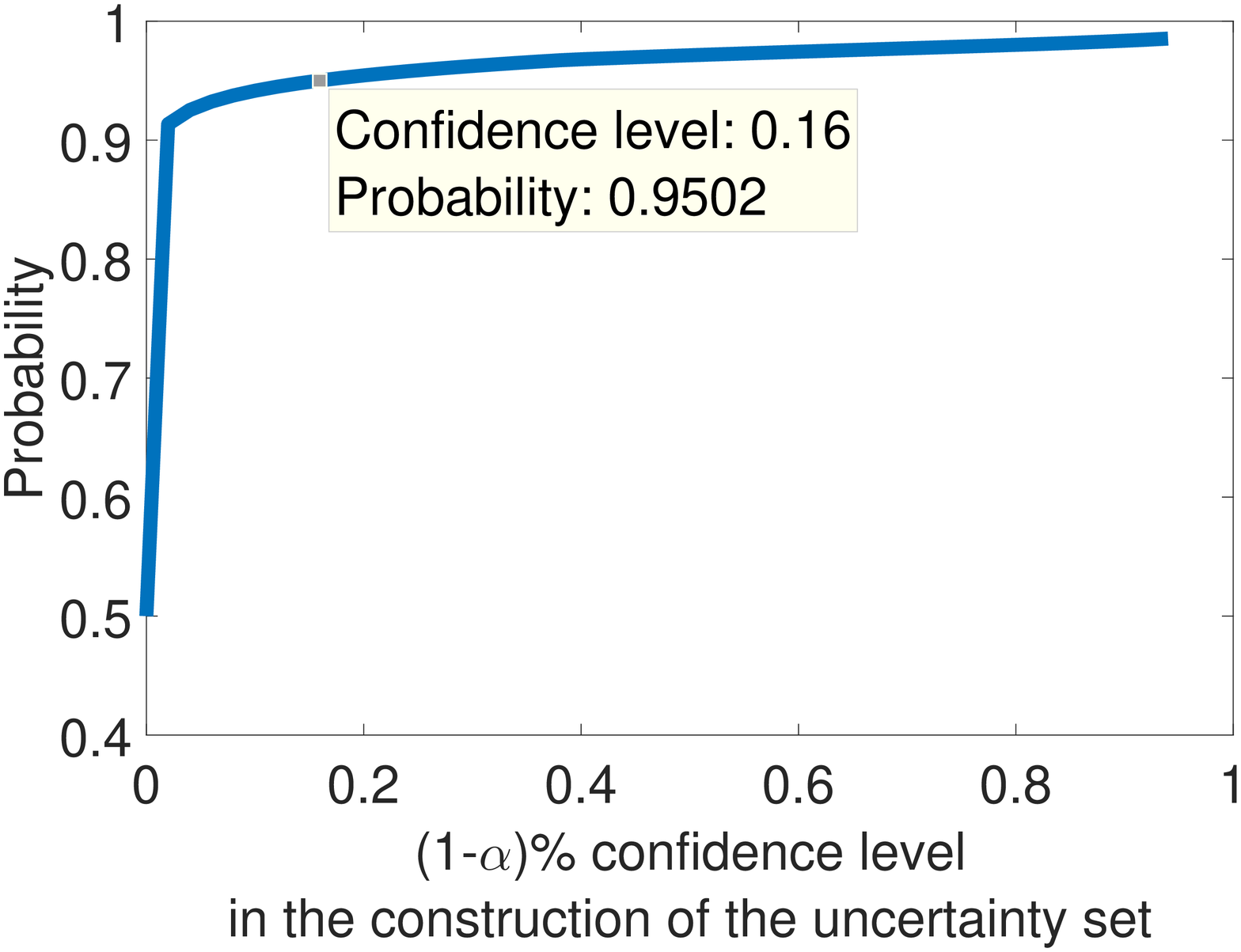}
 		\caption{5 industries}
 	\end{subfigure}
 	\begin{subfigure}{0.5\linewidth}
 		\centering
 		\includegraphics[width=8cm, height=5cm]{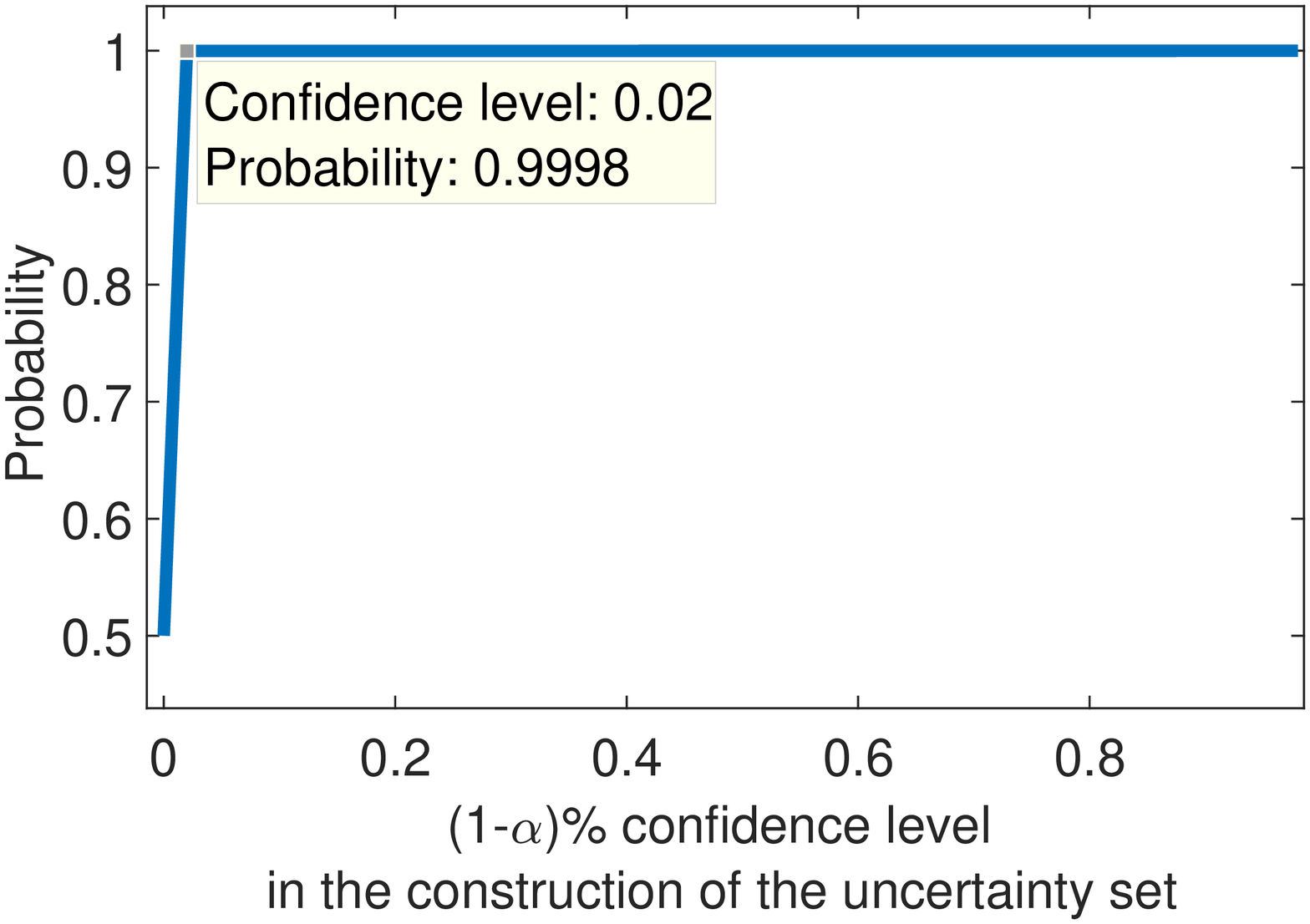}
 		\caption{30 industries}
 	\end{subfigure}
 	\caption{ The horizontal axis presents the value of $(1-\alpha)\%$, the confidence level used in the uncertainty set for data with 5 and 30 industries. The vertical axis presents the value of $Prob(\mu^T\omega^*- 3 \omega^{*^T} \Sigma \omega^*\geq z^*)$, where $\omega^*,z^*$ are the robust solution and corresponding objective value for (\ref{robust risk aversion portfolio}), respectively, with the uncertainty set (\ref{Eq: the new statistical uncertainty set}) and different $\alpha$. The plots are constructed by considering multiplications of $0.02$ in $[0,1]$ as values of $\alpha$.}\label{Fig: 5-30 industries aplha vs z_a for portfolio opt.}
 \end{figure}
 As shown in Figure \ref{Fig: 5-30 industries aplha vs z_a for portfolio opt.}, in order to be sure that the constraint $\mu^T\omega^*- \lambda \omega^{*^T} \Sigma \omega^*\geq z^*$ is satisfied with probability of at least $95\%$, one can reduce the confidence level used in the construction of the uncertainty set from $95\%$ to $16\%$ for the 5 industries case, and to $2\%$ for the 30 industries case.
 
 We emphasize that even though the confidence level of $2\%$ seems to make the uncertainty set much smaller than the one corresponding to $95\%$, this does not happen for the 30 industries case since we have $\sqrt{\frac{\chi^2_{{rank}(V),0.95}}{n}}=0.8723$ and $\sqrt{\frac{\chi^2_{rank(V),0.02}}{n}}=0.7751$, where $n=720$ and $rank(V)=495$.
 \begin{remark}
 	The 2\% confidence level was achieved from Figure \ref{Fig: 5-30 industries aplha vs z_a for portfolio opt.}(b), which was plotted by discretizing $[0,1]$ into the set of points starting from $0$ with the step size of $0.02$. Hence, choosing a smaller step size may result in a smaller confidence level. However, the uncertainty set will not be much smaller than the one constructed by $2\%$ confidence level, as it can be easily checked that even for $10^{-23}\%$ confidence level we have $\sqrt{\frac{\chi^2_{rank(V),10^{-25}}}{n}}=0.5710$, where $n=720$ and $rank(V)=495$.  	
 	\hfill \qed
 \end{remark}
 
 We considered both data sets with 5 and 30 industries in our numerical experiments; however, due to similarity in the results, we present the results of considering only data set with 30 industries. 
 
 After solving the \textit{Nominal problem}, the \textit{Robust problem} considering the uncertainty set with $95\%$ confidence level, the \textit{Robust problem} considering the uncertainty set with $2\%$ confidence level, and the \textit{Chance problem} with $95\%$ confidence level, we compare the solutions in three ways:
 \begin{enumerate}[(i)]
 	\item evaluating the solutions with respect to the nominal values;
 	\item evaluating the solutions with respect to their worst-case scenarios in the uncertainty set constructed with $95\%$ confidence level;
 	\item evaluating the solutions with respect to their worst-case scenarios in the uncertainty set constructed with $2\%$ confidence level.
 \end{enumerate}
 Table  \ref{Table for 30 industries} presents the evaluations of the solutions.
 \begin{table}
 	\begin{center}
 		\begin{tabular}{|c|c|c|c|c|c|}
 			\hline
 			&             &    solution of                        & \multicolumn{2}{c|}{  solution of   } &  solution of\\
 			&             &  \textit{Nominal problem}  & \multicolumn{2}{c|}{\textit{Robust problem} (\ref{robust risk aversion portfolio}) }                               &  \textit{Chance problem} \\
 			&&(\ref{Risk aversion portfolio problem})&\multicolumn{2}{c|}{  with confidence level}&(\ref{Eq: eqivalent risk aversion portfolio with chance constraint})\\
 			&             &                                                     &   $\hspace{0.4cm}95\%\hspace{0.4cm}$ &  \begin{tabular}{c}  $2\%$\end{tabular}                          &  \\ \hline
 			\multirow{3}{*}{Nominal case}  & Obj. value  &                        \textbf{-35.57}                       &  -39.24   &                              -38.84                               &                                      -35.59                                       \\
 			&  \textit{Ann. risk}  &                        \textit{12.03}                        &  \textit{12.63}   &                             \textit{ 12.56}                               &                                      \textit{12.03 }                                      \\
 			& \textit{Ann. return} &                      \textit{7.02}                       & \textit{7.29} &                             \textit{7.28}                             &                                     \textit{6.98}                                     \\ \hline
 			\multirow{3}{*}{$\begin{matrix}
 				\mbox{Worst-case with }  \\
 				\mbox{confidence level } \\
 				95\%
 				\end{matrix}$}   & Obj. value  &                       -77.47                       &    \textbf{-51.11}     &                           -51.12                           &                                      -55.01                                        \\
 			&  Ann. risk  &                        \textit{15.77}                        &     \textit{13.22}     &                            \textit{13.20}                            &                                      \textit{13.42}                                      \\
 			& Ann. return &                     \textit{-183.10}                    &  \textit{-89.14}   &                          \textit{-90.49}                         &                                   \textit{-119.87}                                      \\ \hline
 			\multirow{3}{*}{$\begin{matrix}
 				\mbox{Worst-case with }\\
 				\mbox{confidence level }\\2\%\end{matrix}$}  & Obj. value  &                       -75.00                      &   -50.39  &                             \textbf{-50.38}                             &                                      -53.89                                     \\
 			&  \textit{Ann. risk}  &                        \textit{15.57 }                       &  \textit{13.18}    &                             \textit{13.16}                              &                                      \textit{13.33}                                         \\
 			& \textit{Ann. return} &                     \textit{-172.72}                      & \textit{-83.91} &                           \textit{-85.20}                           &                                    \textit{-113.37}                                      \\ \hline
 		\end{tabular}
 	\end{center}
 	
 	\caption{Comparison among the solutions of the \textit{nominal problem} (\ref{Risk aversion portfolio problem}),  the \textit{Robust problem} (\ref{robust risk aversion portfolio}) considering the uncertainty set with $95\%$ confidence level, the \textit{Robust problem} (\ref{robust risk aversion portfolio}) considering the uncertainty set with $2\%$ confidence level, and (\ref{Eq: eqivalent risk aversion portfolio with chance constraint}) in three way: The first block row with results is the nominal evaluation of the solutions. The second and third block rows are the evaluation of the solutions with respect to their worst-case scenarios in uncertainty sets $95\%$, and $2\%$ confidence level, respectively. The results are by considering the data for 30 industries. 
 		The \textbf{bold numbers} shows the best objective value in each scenario. The annualized return and risk are in \textit{italics} and not individually optimized.}
 	\label{Table for 30 industries}
 \end{table}		
 In the first block row (with results), the evaluation is done using the nominal scenario. The objective value of the \textit{Nominal problem} is the highest. The worst objective value in this row is corresponding to the solution of the \textit{Robust problem} considering the uncertainty set with $95\%$ confidence level. This solution is immunized against more scenarios than the others.
 
 The second block row is the evaluation of the solutions considering their worst-case scenario in the uncertainty set constructed by $95\%$ confidence level. This implies that the solution of the \textit{Robust problem} with this uncertainty set has the highest objective value, because the solution is immunized against all scenarios in the uncertainty set; however,  other solutions are immunized against all scenarios in a subset of the uncertainty set.
 The third block row has the same interpretation, where the scenario is chosen in the uncertainty set with confidence level $2\%$.

 Table \ref{Table for 30 industries} shows that even though all solutions have close annualized returns and risks in the nominal scenario, the solutions of \eqref{robust risk aversion portfolio} have extremely better returns and risks in the included worst-case scenarios. 
 
 
 Proposition \ref{Proposition: chance relaxation is better than robust} states that a solution of \eqref{robust risk aversion portfolio}, denoted by $(\bar{ \omega},\bar{W})$, is more conservative than a solution of \eqref{Eq: eqivalent risk aversion portfolio with chance constraint}, denoted by $(\tilde{ \omega},\tilde{W})$. This means $(\bar{ \omega},\bar{W})$ is safeguarded against more scenarios (all scenarios in the uncertainty set $\bar{ \mathcal{U}}$) than $(\tilde{ \omega},\tilde{W})$. Therefore, as the last column of Table \ref{Table for 30 industries} shows,  the objective values of $(\tilde{ \omega},\tilde{W})$ at their worst-case scenarios in $\bar{ \mathcal{U}}$ are worse than the ones for $(\bar{ \omega},\bar{W})$.
 %
 
 %
 
 \subsection{Least-squares problems with uncertainties}\label{Sec:least square}
 This subsection contains applications of the results of Section \ref{Sec: inner approximation} to two well-known problems, namely a norm approximation and a linear regression problem. 
 \subsubsection{Norm approximation with uncertainty in the coefficients}  \label{Sec: norm approximation}
 The norm approximation $\min_{\adj\in\mathbb{R}^n} \Eucledianorm{A\adj-b}$ tries to find the closest vector to $b\in \mathbb{R}^m$ in the range of the linear function $A\adj$.
 The solution to this problem can be sensitive even to small errors in $A$ or $b$. To detect this, one can 
 analyze the condition number of the matrix $A$ and check the sensitivity of the nominal solution to a perturbation in $A$, see, e.g., Chapter 7 in  \cite{higham2002accuracy}. 
 If the condition number is large, then the solution might be sensitive to a small error in $A$ or $b$, \ahmadreza{hence} not reliable. In this subsection we are using the results of Section \ref{Sec: inner approximation} to deal with this problem. 

 Consider the \ahmadreza{uncertain} norm approximation $\min_\adj \Eucledianorm{(\bar{A}+\Delta)\adj-b}$, where $\Delta \in \uncertset\subseteq\mathbb{R}^{m\times n}$ reflects the uncertainty in $\bar{ A}$. This problem is equivalent to 
 $
 \min_{\adj\ahmadreza{\in\mathbb{R}^n}} \adj^T\left(\bar{A}+\Delta\right)^T\left(\bar{A}+\Delta\right)\adj+2b^T\left(\bar{A}+\Delta\right)\adj+b^Tb.
 $ 
 Now using the results of Section \ref{Sec: inner approximation}, upper and lower bounds on the robust optimal value of this problem are obtained by solving
 \begin{equation}\label{Eq: robust inner norm approximation}
 \min_{\uncertA,\adj}\left\{\trace{( \bar{ A}^T\bar{A}+\Omega^2\identity{n})}{\uncertA}+\delta^*_{\uncertset}(2 \uncertA \bar{ A}^T-2b\adj^T)-2b^T\bar{ A}\adj+\Eucledianorm{b}^2:
 \hspace{0.2cm}\left[ \begin{matrix}\uncertA &\adj\\\adj^T&1\end{matrix}\right]\succeq \zero{n+1\times n+1}\right\},
 \end{equation}
 and 
 \begin{equation}\label{Eq: robust outer norm approximation}
 \min_{\uncertA,\adj}\left\{\trace{ \bar{ A}^T\bar{ A}}{\uncertA}+\delta^*_{\uncertset}(2 \uncertA \bar{ A}^T-2b\adj^T)-2b^T\bar{ A}\adj+\Eucledianorm{b}^2:
 \hspace{0.2cm}\left[ \begin{matrix}\uncertA &\adj\\\adj^T&1\end{matrix}\right]\succeq \zero{n+1\times n+1}\right\},
 \end{equation}
 respectively.

 \ahmadreza{For our numerical experiments,} we construct randomly generated problems with ill-conditioned $\bar{ A}$ as follows: we fix \ahmadreza{$n=100$} and generate randomly a matrix $U\in (0,1)^{n\times n}$ and a vector $ b\in(0,1)^n$. Also, we randomly generate an integer $i$ in $\{1,...,n-1\}$ and construct a diagonal matrix $D$ whose first $i$ diagonal entries are randomly chosen in $(-5,5)$ and the remaining diagonal entries are randomly chosen in $(0,10^{-8})$. Then, we set $\bar{ A}:=U^TDU$. \ahmadreza{Using this procedure, we generate 20 ill-conditioned $\bar{ A}$ matrices with condition numbers in the interval $[10^{15},10^{18}]$. Moreover, we generate uniformly distributed pseudorandom matrices $B^1,B^2\in\{0,1\}^{n\times n}$, and an integer number $K\in\{1,2,...,n^2\}$ using MATLAB built-in function ``randi''. Then,  }we solve the norm approximation problems \ahmadreza{using the generated matrices and the budget-type uncertainty set, proposed in \cite{bertsimas2004price}: 
 	\begin{equation}\label{Eq:uncertainty set with two constraints}
 	\uncertset=\left\{\Delta\in \mathbb{R}^{n\times n}: \;\Vert\Delta\Vert_\infty\leq \rho \right\}\cap \left\{\Delta\in\mathbb{R}^{n\times n}:\enskip \Vert B^k\circ \Delta\Vert_1\leq K\rho, \enskip k=1,2  \right\},
 	\end{equation}	
 	for some $\rho>0$. For this uncertainty set, one can derive $\delta^*_{\uncertset}(U)$ using Lemma \ref{Lemma: composition support}.\eqref{hadamard product}, \ref{Lemma: composition support}.\eqref{intersection}, and \ref{Lemma: Support functions of some uncertainty sets}.\eqref{general norm}. It is worth noting that the constructed uncertain norm approximation problems contain $100\times100=10,000$ uncertain parameters, and hence obtaining an exact optimal robust solution is computationally intractable.  }

 \begin{figure}[t]
 	\begin{subfigure}{0.5 \linewidth}
 		\centering
 		\includegraphics[width=9.0cm, height=5cm]{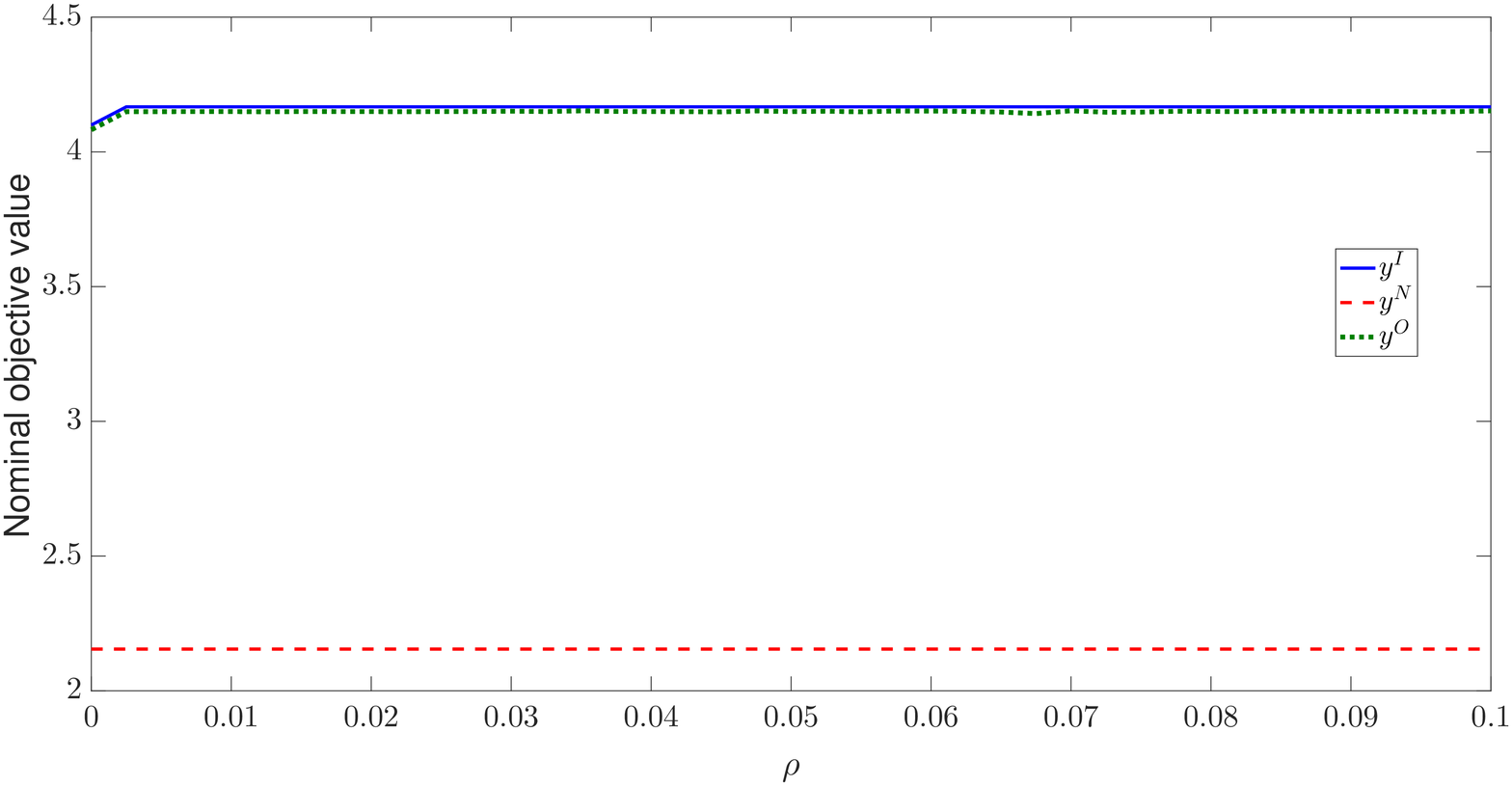}
 		\caption{}\label{Figure: norm approximation nominal average}
 	\end{subfigure}
 	\begin{subfigure}{0.5 \linewidth}
 		\centering
 		\includegraphics[width=9.0cm, height=5.0cm]{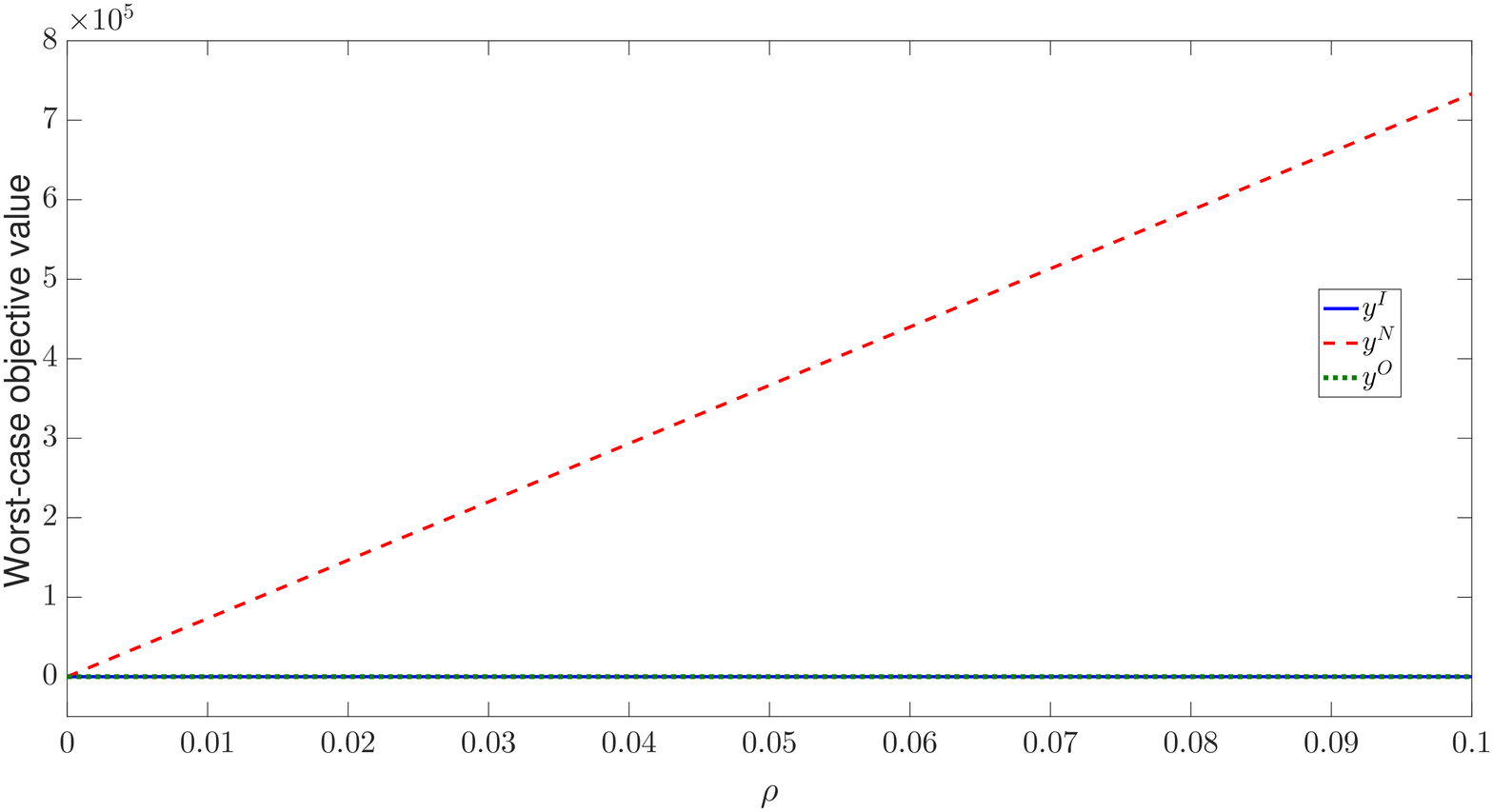}
 		\caption{}\label{Figure: norm approximation worst average}
 	\end{subfigure}
 	\caption{The average behavior \ahmadreza{of the objective values} of $y^I$, $y^O$\ahmadreza{,} and $y^N$ related to \ahmadreza{20} randomly generated norm approximation problems. (a) The nominal objective value is computed by $\Eucledianorm{\bar{ A}\adj-b}$. 
 		(b) The worst-case objective value is computed by $\Eucledianorm{(\bar{A}+\Delta^*)\adj-b}$, where $\Delta^*$ is the worst-case scenario corresponding to $y^N$, $y^I$ or $y^O$. \ahmadreza{Notice that the scales of the vertical axes in (a) and (b) are different.} The solid blue, red dashed, and green \ahmadreza{dotted} curves correspond to $y^I$, $y^N$, and $y^O$, respectively.
 	}\label{Figure: norm approximation average}
 \end{figure}
 
 \ahmadreza{We analyze the performance of the solutions by comparing the objective values of the solutions $y^N$, $y^I$, and $y^O$ for both the nominal matrix $\bar{ A}$ and a worst-case matrix $\bar{ A}+\Delta^*$ corresponding to the vector $y$, where $\Delta^*$ is constructed using the algorithm proposed in Appendix \ref{Appendix:heuristic}.
 }
 
 \ahmadreza{Figure \ref{Figure: norm approximation average} provides a visualization of the average performance of $y^N$, the nominal solution, $y^I$, the solution of \eqref{Eq: robust inner norm approximation}, and $y^O$, the solution of \eqref{Eq: robust outer norm approximation} for different scenarios and values of $\rho\in [0,0.1]$. We check Assumption \eqref{assumption: positive semidefinitness of a part} by solving \eqref{Eq:approximation of DC formulation of assumption psd for all} and find that this assumption holds when $\rho<0.02$ for all instances except one. When this assumption is not satisfied, $y^O$ is just an approximated robust solution, which is more robust than the nominal solution, and \eqref{Eq: robust outer norm approximation} is no longer a lower bound for $$\min_{y\in\mathbb{R}^n}\max_{\Delta\in\mathcal{Z}} \adj^T\left(\bar{A}+\Delta\right)^T\left(\bar{A}+\Delta\right)\adj+2b^T\left(\bar{A}+\Delta\right)\adj+b^Tb.$$
 	One of the important observations from Figure \ref{Figure: norm approximation average} is that even though the constructed matrices are nonsingular, and hence the true nominal objective value is zero, the solver is not able to find the true optimal solution of the nominal problem because of the large condition number of the matrix $\bar{ A}$. 
 	Furthermore, despite the small difference in the average performance of the solutions in the nominal case, the average performance of $y^N$ in its worst-case scenario is extremely worse than the performance of $y^I$ and $y^O$. For instance, for $\alpha=0.7$, the average value of $\Eucledianorm{(\bar{A}+\Delta_N^*)\adj^N-b}$ is $5.13\times 10^{5}$ whereas the average values of $\Eucledianorm{(\bar{A}+\Delta_I^*)\adj^I-b}$ is $4.22$ and the average value of $\Eucledianorm{(\bar{A}+\Delta_O^*)\adj^O-b}$ is $4.29$, where $\Delta_N^*$, $\Delta_I^*$, and $\Delta_O^*$ are the worst-case scenarios corresponding to $y^N$, $y^I$, and $y^O$, respectively. This implies that in average $(\bar{ A}+\Delta_N^*)y^N$ is a point in the range of $\bar{ A}+\Delta_N^*$ that is far from $b$, as $b\in(0,1)^n$.   }
 \begin{figure}[t]
 	\centering
 	\includegraphics[width=14cm, height=12cm]{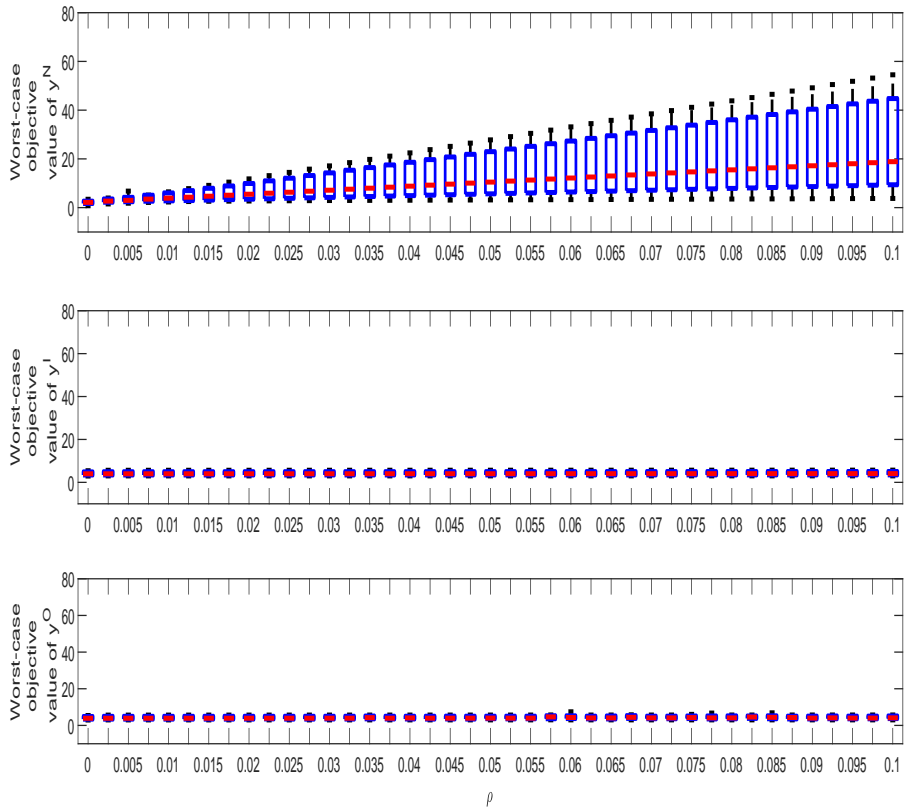}
 	\caption{\ahmadreza{Box plot of the worst-case objective values of $y^N$, $y^I$, and $y^O$ for 20 randomly generated norm approximation problems with $\rho\in[0,0.1]$. The boxes are representing the values between the first and the third quartile. The outliers are left out in the figures to have a better comparison (cf. the main text where more is explained).} }\label{Figure: norm approximation Box plot of worst-case}
 \end{figure}
 
 \ahmadreza{Figure \ref{Figure: norm approximation Box plot of worst-case} provides the box plot of the objective values of the solutions with respect to their worst-case scenarios, where, for each $\rho$, the box represents the values between the first and third quartile, the dashed line above and below each box indicates the range of objective values excluding the outliers, and the red line in each box represents the second quartile. }
 
 \ahmadreza{As this figure shows, the variance of the worst-case objective values of $y^I$ and $y^O$ does not change much as $\rho$ increases. However, the worst-case objective value of $y^N$ significantly changes when $\rho$ increases. This shows the robustness of $y^I$ and $y^O$ against small changes in the components of $\bar{ A}$, whereas $y^N$ is very sensitive to these changes. }
 
 \ahmadreza{   Furthermore, the comparison between Figures \ref{Figure: norm approximation average}.(b) and \ref{Figure: norm approximation Box plot of worst-case} shows that the extremely high average value of $\Eucledianorm{(\bar{A}+\Delta_N^*)\adj^N-b}$ is because of some outliers with extremely high values. However, even after removing the outliers, $y^I$ and  $y^O$ outperform $y^N$.
 }

 \begin{figure}[t]
 	\centering
 	\includegraphics[width=10cm, height=5.3cm]{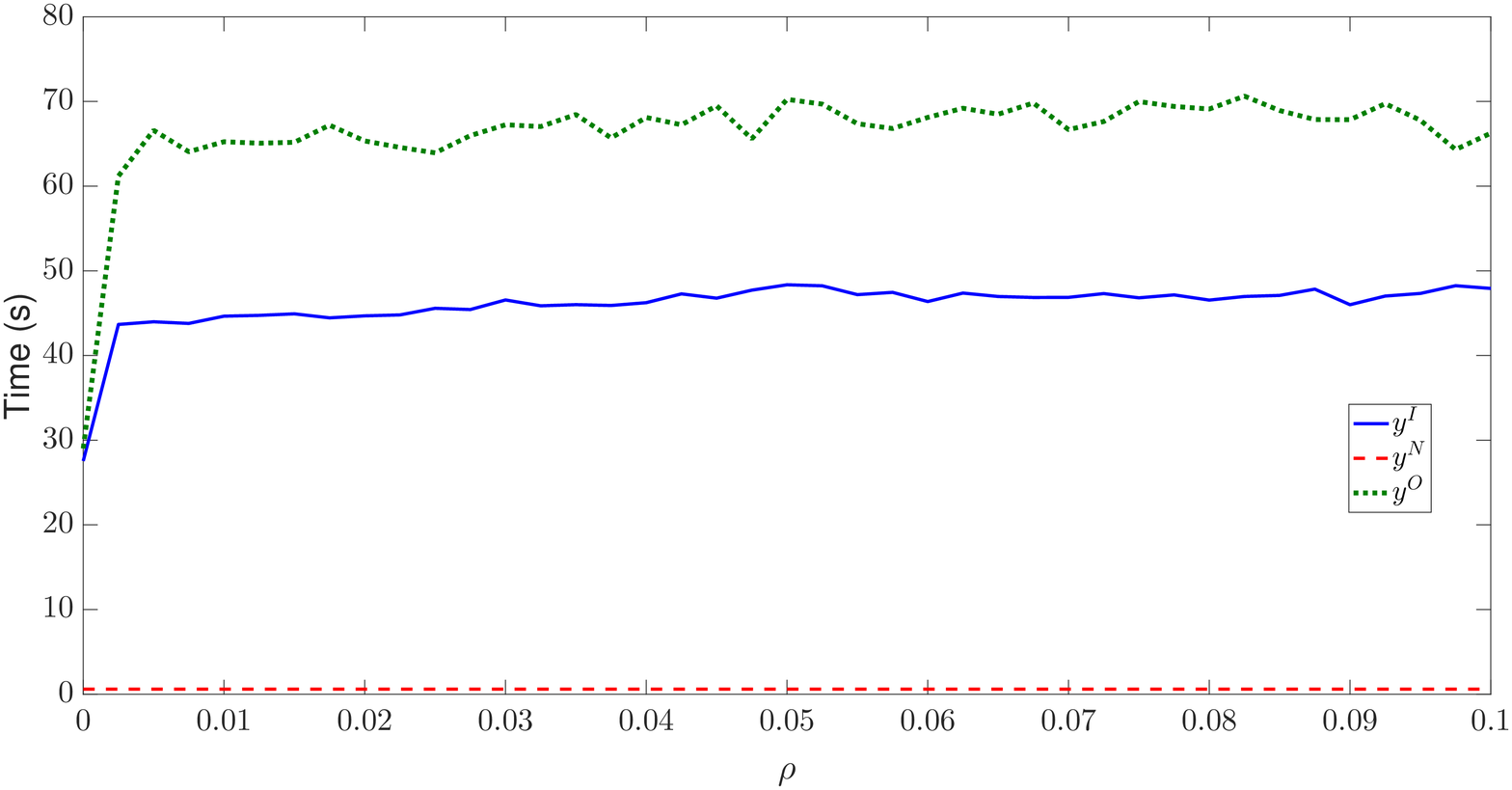}
 	\caption{\ahmadreza{Average time (in seconds) spent by MOSEK to obtain $y^N$, $y^I$, and $y^O$ for 20 randomly generated norm approximation problems. The solid blue, red dashed, and green dotted curves correspond to $y^I$, $y^N$, and $y^O$, respectively. }}\label{Figure: norm approximation time}
 \end{figure}
 \ahmadreza{Figure \ref{Figure: norm approximation time} provides the average time (in seconds) taken by MOSEK to solve the nominal problem as well as \eqref{Eq: robust inner norm approximation} and \eqref{Eq: robust outer norm approximation} to obtain $\adj^N$, $\adj^I$, and $\adj^O$, respectively. We emphasize that even though \eqref{Eq: robust inner norm approximation} and \eqref{Eq: robust outer norm approximation} have only one constraint, which is a linear matrix inequality, we need $O(n^2)$ more variables and constraints to pass the optimization problems to the solver. This is the reason that we see a difference between the time spent to get $y^N$ with the one for $y^I$ and $y^O$ when $\rho=0$.  }

 \subsubsection{Robust linear regression with data inaccuracy}\label{Sec: regression line}
 Another application of the results of this paper is finding a robust linear regression of a dependent variable $Y$ and a vector of independent variables $X$ that are highly collinear. 
 For a data set with $n$ linearly independent variables and $m$ data points, a mathematical formulation of finding the regression line is 
 \begin{equation}\label{Eq: regression line}
 \min_{w,c,b}\left\{\Eucledianorm{w},:
 \hspace{.3cm} w_i=\sum_{j=1}^{n}X_{ij}c_j+b-Y_i,\; \forall i=1,...,m\right\},
 \end{equation}
 where $X_{ij}$ is the $i$-th observed value of the $j$-th independent variable and $Y_i$ is the value of the dependent variable in the $i$-th observation. 
 
 \ahmadreza{For our numerical experiment, we use the dataset proposed in \cite{Candanedo2017datadriven}, which is used to create a regression model of appliances energy consumption in a low-energy house located in Stambruges, Belgium. The dataset consists of 19,735 observations of 26 continuous measurable variables, 10 of which are temperatures of different parts of the house. The description of the variables can be found in Table 2 in \cite{Candanedo2017datadriven}. In this section, we analyze the performance of our results in acquiring a robust linear regression model to predict the appliances energy consumption.}

 To reformulate (\ref{Eq: regression line}) into the form (\ref{concave uncertainty conic quadratic}), let $\bar{ A}\in \mathbb{R}^{m\times(n+2)}$ be a matrix whose collection of the first \ahmadreza{$n=25$} columns is the matrix $X$ consisting of the observations corresponding to all variables \ahmadreza{except the appliances energy consumption, the $(n+1)st=26$th column corresponds to the observations of the appliances energy consumption,} and the components of the last column are all ones. Then, problem (\ref{Eq: regression line}) is equivalent to 
 $
 \min_{\adj\in \mathbb{R}^{(n+2)}}\{\Eucledianorm{\bar{ A}\adj}:\; \adj_{n+1}=-1\}.
 $
 Solving this problem results in the nominal solution $\adj^N$. \ahmadreza{The condition number of $\bar{ A}$ is $1.16\times 10^5$}. This means that the nominal solution might be sensitive to an error in $\bar{ A}$. Let us assume that the maximal inaccuracy in the coefficients of the first $(n+1)$ columns of $\bar{ A}$ is $1\%$, \ahmadreza{and the aggregated inaccuracy in the temperature data cannot exceed $0.1\%$ of the aggregated values. Hence, we consider the following uncertainty set:
 	$$
 	\uncertset=\left\{\Delta\in \mathbb{R}^{m\times (n+2)}:\;
 	\vert\sum_{j\in\mathcal{J}}\sum_{i=1}^m\Delta_{ij}\vert\leq \rho
 	, \;
 	|\Delta_{ij}|\leq \bar{\rho}_j,\;\Delta_{i(n+2)}=0, \begin{matrix}i=1,...,m,\\j=1,...,n+1 \end{matrix} \right\},
 	$$	
 	where $\mathcal{J}$ is the set containing the indices of the temperature columns and $\rho=0.001\sum_{j\in\mathcal{J}}\sum_{i=1}^{m}\vert \bar{ A}_{ij}\vert$, $\bar{\rho}_j=0.01\max_{i=1,...,m}\vert \bar{ A}_{ij}\vert$, $j=1,...,n+1$.}
 To obtain the value of $\Omega$ in Assumption \eqref{assumption: upper bound on the uncertainty}, we first notice that
 \begin{equation}\label{Eq:omega for budgeted set}
 \sup_{\Delta\in \uncertset}\Vert\Delta\Vert_{2,2}\leq \sup_{\begin{tiny}
 	\begin{matrix}
 	\Delta: \\ \vert\Delta_{ij}\vert\leq \bar{\rho}_j\\  \ahmadreza{i=1,...,m\atop j=1,...,n+1.}
 	\end{matrix}
 	\end{tiny}} \sup_{x:\Eucledianorm{x}=1} \Eucledianorm{\Delta x}
 =\sup_{x:\Eucledianorm{x}=1}\sup_{\begin{tiny}
 	\begin{matrix}
 	\Delta: \\ \vert\Delta_{ij}\vert\leq \bar{\rho}_j\\ \ahmadreza{i=1,...,m\atop j=1,...,n+1.}
 	\end{matrix}
 	\end{tiny}}  \Eucledianorm{\Delta x} 
 =\sup_{x:\Eucledianorm{x}=1} \ahmadreza{ \sqrt{m}\Vert[ \bar{\rho}_jx_j]_{j=1,...,n+1}\Vert_1}
 \end{equation}
 \ahmadreza{and due to symmetry we can reformulate  the far right optimization problem in \eqref{Eq:omega for budgeted set} to $$\sup_{x:\Eucledianorm{x}\leq1} \ahmadreza{ \sqrt{m}\sum_{j=1}^{n+1} \bar{\rho}_jx_j}.$$ So, we use 
 	$$
 	\Omega:=\sqrt{m}  \begin{aligned}\max_{x\in\mathbb{R}^{n+2}:\Vert x\Vert_2\leq1}\hs\sum_{j=1}^{n+1}\bar{ \rho}_jx_j.
 	\end{aligned}
 	$$}
 For this instance, Assumption (\ref{assumption: positive semidefinitness of a part}) does not hold. Therefore, we only consider the inner approximation 
 $$
 \min_{\begin{tiny}
 	\begin{matrix}
 	\uncertA\in \mathbb{R}^{(n+2)\times (n+2)}\\\adj\in \mathbb{R}^{(n+2)}\\\scalervar\in \mathbb{R}
 	\end{matrix}
 	\end{tiny}}\left\{\trace{( \bar{ A}^T\bar{ A}+\Omega^2\identity{n+2})}{\uncertA}+\delta^*_{\uncertset}(2\bar{ A} \uncertA )+\frac{\scalervar}{4}:
 \hspace{0.2cm}\left[ \begin{matrix}\uncertA &\adj\\\adj^T&\scalervar\end{matrix}\right]\succeq \zero{n+1\times n+1},\hs y_{n+1}=-1\right\},
 $$ 
 to obtain a robust solution $y^I$. 
 \ahmadreza{We use $75\%$ of the observations in the dataset to construct the regression model and the remaining ones to test the performance of the regression lines. To obtain the robust regression line, denoted by $y^I$, we solve the inner approximation, which takes around $276.5$ seconds to be solved. Moreover, the nominal regression line, denoted by $y^N$, is obtained in around $1.1$ seconds by solving \eqref{Eq: regression line}. }
 \begin{figure}
 	\centering
 	\includegraphics[width=18cm, height=7cm]{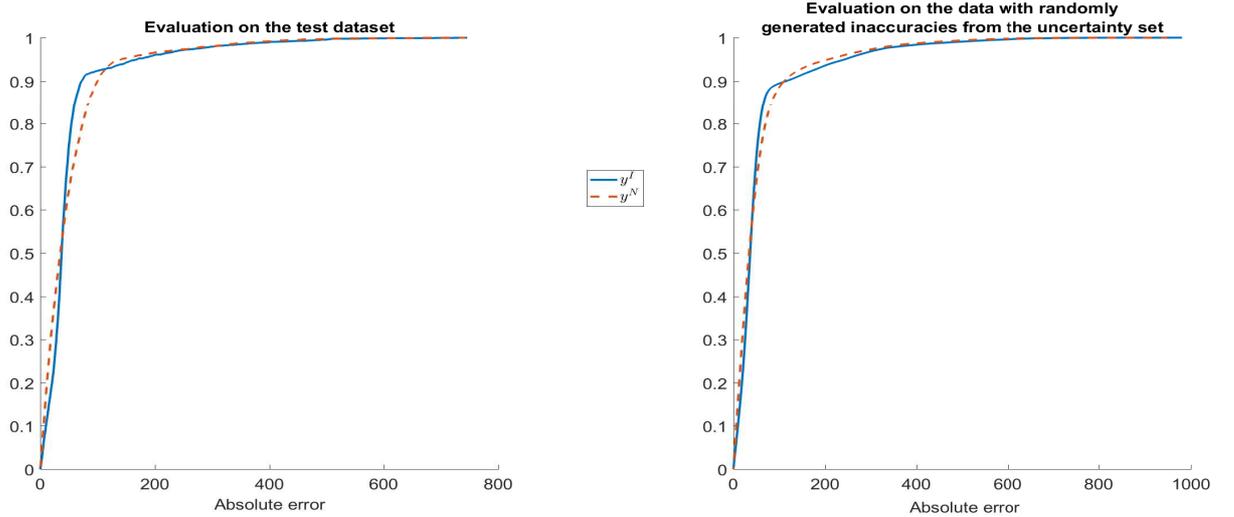}
 	\caption{\ahmadreza{The empirical cumulative distribution functions (ECDF) of the absolute errors of two linear regressions: the solid blue line corresponds to $y^I$ obtained from our inner approximation problem and the dashed red line corresponds to $y^{N}$ obtained by solving \eqref{Eq: regression line}. 		 The left and right figures illustrate the ECDF of the absolute errors of the two linear regressions on the test dataset and the randomly generated scenarios, respectively.}  }\label{Figure: Regression line distribution test }
 \end{figure}
 \ahmadreza{ Figure \ref{Figure: Regression line distribution test } shows the empirical cumulative distribution function (ECDF) of the absolute errors of each regression lines on: (i) the rest of the 25\% observations in the dataset, called the test dataset, and (ii) on 100 randomly generated scenarios $\Delta$ from the uncertainty set $\mathcal{Z}$. We use each random scenario to construct a possible inaccuracy noise that is neglected in the data.  As one can see in Figure \ref{Figure: Regression line distribution test }, the performances of $y^I$ and $y^N$ are close while measured on the test dataset and randomly generated scenarios. The main distinction of the performances are in the skewness of the empirical distributions of the absolute errors. As can be seen, both distributions are right-skewed but the distribution of the absolute errors of the regression line $y^N$ has a higher skewness. 
 }	 

 \begin{table}[h]
 	\ahmadreza{	\centering
 		\begin{tabular}{|c|c|c|}
 			\hline
 			&  $y^I$   &  $y^N$   \\ \hline
 			$\Eucledianorm{\left(\bar{ A}+\Delta^*_{I}\right)y}$ & 13,714.6 & 15,104.0 \\ \hline
 			$\Eucledianorm{\left(\bar{ A}+\Delta^*_{N}\right)y}$ & 13,559.7 & 15,272.3 \\ \hline
 		\end{tabular}
 		\caption{\ahmadreza{Performances of the inner approximation solution, $y^I$, and the nominal solution, $y^{N}$, on the worst-case scenarios. $\Delta^*_{I}$ and $\Delta^*_{N}$ are the scenarios generated in the uncertainty set using the algorithm in Appendix \ref{Appendix:heuristic} for $y^I$ and $y^{N}$, respectively.} }\label{Table:comparison of legression models}}
 \end{table}

 \ahmadreza{Table \ref{Table:comparison of legression models} shows the comparison between the worst-case performances of $y^I$ and $y^{N},$ where $\Delta^*_{I}$ and $\Delta^*_{N}$ are the worst-case scenarios for $y^I$ and $y^{N}$, respectively, obtained using the heuristic algorithm proposed in Appendix \ref{Appendix:heuristic}. As one can see, $y^I$ has a better performance with respect to both worst-case scenarios $\Delta^*_{I}$ and $\Delta^*_{N}$. More specifically, the regression line obtained from the inner approximation results in $10\%$ improvement on the error realized in the worst-case scenarios. This shows the superiority of the robust regression line since the difference between its performances in the randomly generated inaccuracies as well as the test dataset are close to the ones delivered using the nominal solution. }

 \section*{Acknowledgment}
 The research of the first author is partially supported by EU Marie Curie Initial Training Network number 316647 (``Mixed Integer Nonlinear Optimization (MINO)"), \ahmadreza{as well as by the 4TU strategic research and capacity building programme DeSIRE (Designing Systems for Informed Resilience Engineering), as part of the 4TU-programme High Tech for a Sustainable Future (HTSF).}
 \\
 Furthermore, the authors would like to thank the anonymous reviewers and the associate editor for their helpful and constructive comments that greatly contributed to improving the paper.
 
 \begin{appendices}
 	The appendix of this paper consists of \ahmadreza{four} parts. 
 	Part \ref{Sec: appendix proofs} contains the proofs of \ahmadreza{several} lemmas and propositions. In the second part, we provide two simple examples to illustrate the results of Section \ref{Sec:exact formulation}. \ahmadreza{We provide different methods in Appendix \ref{APP: some sets for assumption psd for all} to check Assumptions \eqref{assumption: upper bound on the uncertainty}  and \eqref{assumption: positive semidefinitness of a part}. Finally, we propose a heuristic algorithm in Appendix \ref{Appendix:heuristic} to find a worst-case scenario  in the uncertainty set defined in \eqref{Eq:uncertainty set with two constraints} corresponding to a solution $y$.}

 	\section{Proofs}\label{Sec: appendix proofs}

 	\subsection{Proof of Lemma \ref{Lemma: composition support}}\label{Appendix: proof of the lemma of the support function with some composition rule }
 	
 	(\ref{vec })
 	$
 	\delta^*_{\uncertset}(U)=\sup_{\Delta\in \mathbb{R}^{n \times n}} \left\{\trace{\Delta}{ U^T}:\hspace{.2cm}vec(\Delta)\in \mathcal{U}\right\}
 	=\sup_{\Delta\in \mathbb{R}^{n \times n}}\left\{vec(U)^Tvec(\Delta):\hspace{.2cm}vec(\Delta)\in \mathcal{U}\right\}=\delta^*_{\mathcal{U}}(vec(U)).
 	$
 	\\
 	(\ref{summation})	
 	$
 	\delta_{\uncertset}^*(U)=\sup_{\Delta\in \uncertset}\left\{\trace{\Delta}{ U^T} \right\}=\sup_{\uncert\in \mathcal{U}}\left\{\sum_{i=1}^k\trace{\uncert_i\Delta^i}{ U^T} \right\}
 	$\\$	=\sup_{\uncert\in \mathcal{U}}\left\{\uncert^T\left[\trace{\Delta^i}{ U^T}\right]_{i=1,...,k} \right\}
 	= \delta^*_\mathcal{U}\left(\left[\trace{\Delta^i}{ U^T}\right]_{i=1,...,k}\right).
 	$
 	\\
 	(\ref{general format of affinity})
 	$
 	\delta^*_\uncertset(U)=\sup_{\Delta\in \uncertset}\left\{\trace{\Delta}{U^T} \right\}=\sup_{\Theta\in \mathcal{U}}\left\{\trace{L\Theta R}{U^T} \right\}
 	=\sup_{\Theta\in \mathcal{U}}\left\{\trace{\Theta }{RU^TL} \right\}=\delta^*_\mathcal{U}(L^TUR^T).
 	$\\
 	\eqref{hadamard product}
 	\ahmadreza{
 		Let $\bar{\Delta}\in\mathbb{R}^{n\times n}$ be such that $L\circ \bar{\Delta}\in\mathcal{U}$. If there exist $\bar{ i},\bar{ j}=1,...,n$ for which $L_{\bar{i}\bar{j}}=0$ and $U_{\bar{i}\bar{j}}\neq0,$ then for any integer $k$, set
 		$$\Delta^{k}_{ij}:=\left\{\begin{matrix}\bar{\Delta}_{ij}&\enskip \mbox{if }i\neq\bar{i} \mbox{ or} \enskip j\neq\bar{j},\\
 		ksgn(U_{\bar{i}\bar{j}})\bar{\Delta}_{\bar{i}\bar{j}}&\mbox{ otherwise},
 		\end{matrix}\right.$$ 
 		where $sgn(.)$ is the sign function. This constructs a sequence of matrices $\{\Delta^k\}_{k\in\mathbb{N}}$ for which $L\circ\Delta^k\in\mathcal{U}$ and $\trace{\Delta^k}{U^T}$ goes to $+\infty$ when $k$ tends to $+\infty$.
 		\\
 		Now, assume $U_{ij}=0$ if $L_{ij}=0$, for any $i,j=1,...,n$. Then, 
 		\begin{eqnarray}
 		\delta^*_\uncertset(U)&=&\sup_{L\circ\Delta\in\mathcal{U}}\sum_{i,j=1}^n\Delta_{ij}U_{ij}=\sup_{L\circ\Delta\in\mathcal{U}}\sum_{i,j=1,...,n:\atop L_{ij}\neq0}\Delta_{ij}U_{ij}=\sup_{L\circ\Delta\in\mathcal{U}}\sum_{i,j=1,...,n:\atop L_{ij}\neq0}L_{ij}\Delta_{ij}U_{ij}L^\dagger_{ij}\nonumber\\
 		&=&\sup_{L\circ\Delta\in\mathcal{U}}\trace{L\circ\Delta}{\left(U\circ L^\dagger\right)^T}=\delta^*_{\mathcal{U}}\left(U\circ L^\dagger\right).\nonumber
 		\end{eqnarray}		
 	}
 	(\ref{minkovski sum})$\hs
 	\begin{aligned}
 	\delta^*_\uncertset(U)=\sup_{\Delta \in \uncertset}\trace{\Delta}{U^T}=\sup_{\Delta_i \in \uncertset^i\atop i=1,...,k}\sum_{i=1}^k\trace{\Delta^i}{ U^T}
 	=\sum_{i=1}^k\sup_{\Delta^i \in \uncertset_i}\trace{\Delta^i}{ U^T}=\sum_{i=1}^k\delta^*_{\uncertset_i}(U).
 	\end{aligned}
 	$\\
 	(\ref{intersection}) Similar to the proof of Lemma 9 in \cite{FenchelDuality}.\\
 	(\ref{cartesian product dependent})
 	$
 	\hs\delta^*_\uncertset\left((U_1,...,U_k)\right)=\sup_{\Delta\in \uncertset} \trace{\Delta}{(U_1,...,U_k)^T}=\sup_{\Delta_i\in \uncertset_i\atop i=1,...,k}\trace{(\Delta_1,...,\Delta_k)}{(U_1,...,U_k)^T}
 	=\sup_{\Delta_i\in \uncertset_i\atop i=1,...,k}\trace{\sum_{i=1}^k\Delta_iU_i^T}{}=\sum_{i=1}^k\sup_{\Delta_i\in \uncertset_i}\trace{\Delta_iU_i^T}{}=\sum_{i=1}^k\delta^*_{\uncertset_i}(U_i).
 	$\\
 	(\ref{convex hull})
 	$
 	\delta^*_\uncertset(U)=\sup_{\Delta\in \uncertset}\trace{\Delta}{U^T}=\sup_{\begin{tiny}\begin{matrix}
 		\Delta^i\in \uncertset_i\\
 		\lambda_i\geq 0\\
 		i=1,...,k
 		\end{matrix}
 		\end{tiny}}\left\{\sum_{i=1}^k\lambda_i\trace{\Delta^i}{U^T}:\sum_{i=1}^{k}\lambda_i=1
 	\right\}	=\\ \max_{i=1,...,k}\sup_{\Delta^i\in \uncertset_i}\trace{\Delta^i}{U^T}=\max_{i=1,...,k}\delta^*_{\uncertset_i}(U).
 	$\hfill \qed
 	\subsection{Proof of Lemma \ref{Lemma: Support functions of some uncertainty sets}(\ref{PSD bounded set})}\label{Appendix proof of lemma support function of natural sets}
 	The assumptions imply that 
 	\begin{eqnarray}
 	\delta_\uncertset ^*(U)&=&
 	\sup_{\Delta} \left\{ \trace{\Delta}{U^T}:
 	\; \Delta^l\preceq \Delta \preceq \Delta^u\right\}=\max_{\Delta} \left\{ \trace{\frac{U+U^T}{2}}{\Delta}:
 	\; \Delta^l\preceq \Delta \preceq \Delta^u\right\}\nonumber \\
 	&=&\min_{\Lambda_1,\Lambda_2} \left\{\trace{\Delta^u}{\Lambda_2}-\trace{\Delta^l}{\Lambda_1}:\;\Lambda_2-\Lambda_1=\frac{U+U^T}{2},\; \Lambda_1,\Lambda_2\succeq\zero{n\times n}\right\},\nonumber 
 	\end{eqnarray}
 	where the last equality holds because of conic duality (both problems are strictly feasible).\hfill\qed
 	\subsection{Proof of Lemma \ref{Lemma: SDP matrix norm}$(ii)$}\label{Appendix lemma SDP refomulation of trace norm}	
 	Lemma \ref{Lemma: dual norm}(c) implies that $\Vert U\Vert_{2,2}^2$ is the largest eigenvalue of $UU^T$. Hence, $\Vert U\Vert_{2,2}^2\leq \rho^2$ can be reformulated as $UU^T\preceq \rho^2\identity{n}$, which by using Schur Complement Lemma (see, e.g., Appendix A.5.5 in  \cite{boyd2004convex})  is equivalent  to 
 	$\left[\begin{matrix}
 	\rho^2\identity{n}&U\\ U^T&\identity{n}
 	\end{matrix}\right]\succeq\zero{2n\times 2n}.	$\hfill\qed
 	\subsection{Proof of the statement in Example \ref{Ex: convex quadratic with linear vector uncertainty}}\label{Appendix proof of example with nonegative uncertainty and show no LMI}
 	$\adj\in \mathbb{R}^n$ satisfies (\ref{Eq: example with vector uncertainty}) if and only if 
 	\begin{equation}\label{Eq: lemma with quadratic and nonnegative uncertainty with sup}
 	y^T\bar{ A}y+\sup_{\uncert\in \uncertset}\left\{\uncert^T \left[\adj^TA^i\adj+b^{i^T}\adj\right]_{i=1,...,t}\right\}+\bar{ b}^Ty+c\leq 0.
 	\end{equation}
 	Now, we show that $y\in \mathbb{R}^n$ satisfies (\ref{Eq: lemma with quadratic and nonnegative uncertainty with sup}) if and only if there exists $v\in \mathbb{R}^t$ such that
 	\begin{equation}\label{Eq: lemma with quadratic and nonnegative uncertainty with sup and v}
 	y^T \bar{ A}y+\sup_{\uncert\in \uncertset}\left\{\uncert^T v\right\}+\bar{ b}^Ty+c\leq 0, \;v\geq \left[\adj^TA^i\adj+b^{i^T}\adj\right]_{i=1,...,t}.
 	\end{equation}
 	It is clear that if $\adj \in \mathbb{R}^n$ and $v\in \mathbb{R}^t$ satisfy (\ref{Eq: lemma with quadratic and nonnegative uncertainty with sup and v}) then due to nonnegativity of $\uncert \in \uncertset$, 
 	$$\uncert^Tv\geq \uncert^T\left[\adj^TA^i\adj+b^{i^T}\adj\right]_{i=1,...,t},$$
 	which implies $y\in \mathbb{R}^n$ satisfies (\ref{Eq: lemma with quadratic and nonnegative uncertainty with sup}). Now let $y\in \mathbb{R}^n$ satisfies (\ref{Eq: lemma with quadratic and nonnegative uncertainty with sup}). Then setting $v= \left[\adj^TA^i\adj+b^{i^T}\adj\right]_{i=1,...,t}$ implies $\adj\in \mathbb{R}^n$ and $v\in \mathbb{R}^t$ satisfy (\ref{Eq: lemma with quadratic and nonnegative uncertainty with sup and v}), which can be reformulated as (\ref{Eq: lemmavector uncertainty for quadratic without LMI}).\hfill \qed	

 	\subsection{Proof of Proposition \ref{Proposition: chance relaxation is better than robust}} \label{Proof: proposition chance relaxation is better than robust}
 	Let $\left(\adj, W\right)$ satisfy
 	$$
 	\begin{aligned}
 	&\hat{\mu}^T\adj+\trace{\hat{\Sigma}}{\uncertA}+\sqrt{\frac{\chi^2_{rank(V),1-\alpha}}{n}}\Eucledianorm{\left(\Psi R^{-1}\right)^T\left(\adj \atop svec(\uncertA)\right)}+c\leq 0,\;
 	&\left[\begin{matrix}
 	\uncertA & \adj\\
 	\adj^T&1
 	\end{matrix}\right]\succeq \zero{n+1\times n+1}.
 	\end{aligned}
 	$$
 	We know that $\sqrt{\chi^2_{d,1-\alpha}}\geq z_{1-\alpha},$ for any $d$. Also, we have 
 	$$\trace{W}{\hat{\Sigma}}\geq \trace{\adj \adj^T}{\hat{\Sigma}}=\adj^T\hat{\Sigma}\adj,$$ 
 	where the inequality is because $\hat{\Sigma}\succeq 0_{n\times n}$ and $W\succeq \adj \adj^T$ using Schur Complement Lemma. Therefore, $(\adj,W)$ satisfies 
 	$$
 	\begin{aligned}
 	&\hat{\mu}^T\adj+\adj^T\hat{\Sigma}\adj+\frac{z_{1-\alpha}}{\sqrt{n}}\Eucledianorm{\left(\Psi R^{-1}\right)^T\left(\adj \atop svec(\uncertA)\right)}+c\leq 0,\;
 	&\left[\begin{matrix}
 	\uncertA & \adj\\
 	\adj^T&1
 	\end{matrix}\right]\succeq \zero{n+1\times n+1},
 	\end{aligned}
 	$$
 	which is the same as \eqref{Eq: relaxed chance constraint reformulation}, since $\Psi$ is diagonal.\hfill \qed

 	\section{Some illustrative examples}\label{App: illustrative examples}
 	\begin{example}\label{Ex: Ferbnius norm}
 		Let $\uncertset=\{\Delta\in \mathbb{R}^{n \times n}:\; \Vert\Delta\Vert_F\leq 1\}$ and let the assumptions of Theorem \ref{Th: RC for conic and convex quadratic} hold. Then,  using Theorem \ref{Th: RC for conic and convex quadratic}(\ref{Th: item, convex quadratic}), Lemma \ref{Lemma: dual norm}, and Lemma \ref{Lemma: Support functions of some uncertainty sets}(\ref{general norm}), $\adj$ satisfies (\ref{concave uncertainty convex quadratic}) if and only if there exists $\uncertA\in \mathbb{R}^{n\times n}$ such that
 		$
 		\begin{matrix}
 		\trace{\bar{A}}{\uncertA}+\bar{ b}^T\adj+ c+\left\Vert  \uncertA +  \adj a^T\right\Vert_F\leq 0& \mbox{and }
 		\left[\begin{matrix}
 		\uncertA& \adj\\
 		\adj^T& 1
 		\end{matrix}\right]\succeq\zero{n+1\times n+1}.\\
 		\end{matrix} 
 		$ \hfill \qed
 	\end{example}
 	In the next example, we derive a tractable reformulation of the RC in the form (\ref{concave uncertainty conic quadratic}) with the uncertainty set similar to the one proposed by \cite{delage2010distributionally}.
 	\begin{example}\label{Ex: delage uncertainty}
 		Consider the constraint
 		\begin{equation}\label{Eq: example with delage}
 		\sqrt{\adj^T A(\Delta)\adj}+b(\zeta)^T\adj+c\leq 0\hs \forall(\zeta,\Delta)\in \uncertset,
 		\end{equation}
 		where	$\zeta\in \mathbb{R}^n$, $\Delta\in \mathbb{R}^{n\times n}$ are uncertain parameters, $A(\Delta)=\bar{A}+\Delta $, $b(\zeta)=\bar{b}+D\zeta$,  $\bar{ A},D\in \mathbb{R}^{n\times n}$, $\bar{ b}\in \mathbb{R}^n$, and $\uncertset=\uncertset_1\cap\uncertset_2$,
 		$$
 		\uncertset_1=\left\{(\zeta,\Delta):\left[\begin{matrix}
 		1&\zeta^T\\
 		\zeta&\Delta
 		\end{matrix}\right]\succeq\zero{n+1\times n+1}\right\},\hs \uncertset_2=\left\{(\zeta,\Delta): \Delta^l\preceq \Delta \preceq\Delta^u \right\},
 		$$
 		with given $\Delta^l$ and $\Delta^u$ such that $\Delta^u-\Delta^l\succ\zero{n \times n}$. Also, assume that the assumptions of Theorem \ref{Th: RC for conic and convex quadratic} hold. By Lemma \ref{Lemma: composition support}(\ref{intersection}), 
 		$$
 		\delta^*_\uncertset(U,v)=\min_{U^1,U^2\in \mathbb{R}^{n \times n}\atop v^1,v^2\in \mathbb{R}^n}\left\{\delta^*_{\uncertset_1}(U^1,v^1)+\delta^*_{\uncertset_2}(U^2,v^2):\hs U^1+U^2=U,\;v^1+v^2=v\right\}.
 		$$
 		Following a similar line of reasoning as in the proof of Theorem \ref{Th: RC for conic and convex quadratic}(\ref{Th: item, conic quadratic}), $\adj\in \mathbb{R}^n$ satisfies (\ref{Eq: example with delage}) if and only if there exist   $\uncertA,U^1,U^2 \in \mathbb{R}^{n \times n}$, $v^1, v^2 \in \mathbb{R}^n$ and $\scalervar \in \mathbb{R}$ such that
 		\begin{equation}\label{Eq: example with delage reform. with support}
 		\left\{\begin{matrix}
 		\trace{\bar{A}}{\uncertA}+\bar{ b}^T\adj+ c+\delta^*_{\uncertset_1}(v^1,U^1)+\delta^*_{\uncertset_2}(v^2,U^2)+\frac{\scalervar}{4}\leq 0,\\
 		\left[\begin{matrix} \uncertA&\adj\\\adj^T&\scalervar
 		\end{matrix}\right]\succeq \zero{n+1\times n+1},\;U^1+U^2=\uncertA,\;v^1+v^2=D^T \adj.
 		\end{matrix}\right.
 		\end{equation}
 		Using Lemma \ref{Lemma: Support functions of some uncertainty sets}(\ref{PSD bounded set}), (\ref{Eq: example with delage reform. with support}) is equivalent to
 		$$
 		\left\{\begin{aligned}
 		&\trace{\bar{A}}{\uncertA}+\bar{ b}^T\adj+ c+\trace{\Delta^u}{\Lambda_2}-\trace{\Delta^l}{\Lambda_1}+\frac{\scalervar}{4}+\gamma\leq 0,\\&
 		\;U^1+U^2= \uncertA,\hs
 		\Lambda_2-\Lambda_1=\frac{U^2+U^{2^T}}{2},\hs \Lambda_1,\Lambda_2\succeq \zero{n\times n},	\\&
 		\left[\begin{matrix} \uncertA&\adj\\\adj^T&\scalervar
 		\end{matrix}\right]\succeq \zero{n+1\times n+1},\hs\left[\begin{matrix}
 		\frac{U^1+U^{1^T}}{2}&\frac{1}{2}D^T\adj\\\frac{1}{2}\adj^TD&-\gamma
 		\end{matrix}\right]\preceq\zero{n+1\times n+1},
 		\end{aligned}\right.
 		$$
 		for some $\Lambda_1,\Lambda_2 \in S_n$ and $\gamma \in \mathbb{R}$.	\hfill \qed
 	\end{example}
 	
 	\section{\ahmadreza{How to check Assumptions \eqref{assumption: upper bound on the uncertainty}  and \eqref{assumption: positive semidefinitness of a part}}}\label{APP: some sets for assumption psd for all}
 	
 	\ahmadreza{In this section, we provide methods that can be used to check Assumptions \eqref{assumption: upper bound on the uncertainty}  and \eqref{assumption: positive semidefinitness of a part}, each in a separate subsection. }

 	\subsection{\ahmadreza{Finding $\Omega$ for which Assumption (\ref{assumption: upper bound on the uncertainty}) holds}}
 	Assumption (\ref{assumption: upper bound on the uncertainty}) states that there exists an upper bound $\Omega$ for $\sup_{\Delta \in \uncertset} \Vert \Delta \Vert_{2,2}$. \ahmadreza{This assumption is equivalent to the boundedness of the uncertainty set $\uncertset$. So, checking this assumption can be done easily; however, in our inner approximation, we use the value of $\Omega$. Thus, in this section we provide methods to obtain it.} Notice that $\Vert.\Vert_{2,2}$ is a convex function and the maximization of a convex function over a set, in general, is \textit{NP-hard}. However, \ahmadreza{in Proposition \ref{proposition: Omega for box uncertainty}, we show} how to compute $\Omega$ for the box uncertainty set.

 	In the cases for which $\sup_{\Delta \in \uncertset} \Vert \Delta \Vert_{2,2}$ cannot be computed efficiently, one may use an upper bound for $\Vert.\Vert_{2,2}$ to calculate $\Omega$. \ahmadreza{For instance, one can approximate $\mathcal{Z}$ with the union of simplices or boxes (see, e.g., \cite{ben2016tractable} and \cite{BEMPORAD2004151}) and then find the maximum of $\Vert \Delta \Vert_{2,2}$ over the union by means of Proposition \ref{proposition: Omega for box uncertainty} (for boxes) or by checking the vertices (for simplices).} 
 	
 	\subsection{\ahmadreza{Checking Assumption \eqref{assumption: positive semidefinitness of a part}}}
 	Regarding Assumption \eqref{assumption: positive semidefinitness of a part}, it is mentioned  in Section 8.2 in  \cite{bental2009robust} that finding a robust solution to an uncertain linear matrix inequality, in general, is \textit{NP-hard}. Hence, there is no efficient way, in general,  to check Assumption (\ref{assumption: positive semidefinitness of a part}) exactly. 
 	%
 	%
 	However, there is much research that provides different methods to check Assumption \eqref{assumption: positive semidefinitness of a part}. We refer the reader to the papers  \cite{el1998robust}, \cite{Bertsimas2006uncertainConicConstraints}, and Chapters 8 and 9 of the book  \cite{bental2009robust}. Moreover, the problem may have specific characteristics from which this assumption can be certified. We refer the reader to the similar discussion in Section \ref{Sec:exact formulation}. In the following proposition, \ahmadreza{we provide an equivalent statement to Assumption (\ref{assumption: positive semidefinitness of a part}) for uncertainty sets defined by matrix norms}.
 	
 	\begin{proposition}\label{Proposition: psd to sigma norm}
 		\ahmadreza{Let $
 			\uncertset=\left\{\Delta\in\mathbb{R}^{n\times n}:\Vert\Delta\Vert\leq \rho \right\},
 			$
 			for a general matrix norm $\Vert.\Vert$. Then Assumption (\ref{assumption: positive semidefinitness of a part}) holds if and only if $\sup_{y\in\mathbb{R}^n}\left\{-y^T\bar{ A}^T\bar{ A}y+2\rho\Vert\bar{ A}yy^T\Vert^*\right\}=0$.}
 	\end{proposition}
 	\begin{proof}{Proof.}

 		\ahmadreza{
 			Assumption \eqref{assumption: positive semidefinitness of a part} holds if and only if 
 			\begin{eqnarray}
 			&\forall y\in\mathbb{R}^n\hs \forall\Delta\in\uncertset&\hspace{1cm} y^T\bar{ A}^T\bar{ A}y+2\trace{yy^T\bar{ A}^T}{\Delta}\geq0\nonumber\\
 			\Leftrightarrow	&\forall y\in\mathbb{R}^n&\hspace{1cm} -y^T\bar{ A}^T\bar{ A}y+\sup_{\Delta\in\uncertset}\left\{-2\trace{yy^T\bar{ A}^T}{\Delta}\right\}\leq0\nonumber\\
 			\Leftrightarrow	&\forall y\in\mathbb{R}^n&\hspace{1cm}-y^T\bar{ A}^T\bar{ A}y+\delta_{\uncertset}^*(-2\bar{ A}yy^T)\leq0\nonumber\\
 			\Leftrightarrow	&\forall y\in\mathbb{R}^n&\hspace{1cm}-y^T\bar{ A}^T\bar{ A}y+2\rho\Vert\bar{ A}yy^T\Vert^*\leq0\nonumber\\
 			\Leftrightarrow	&&\sup_{y\in\mathbb{R}^n}\left\{-y^T\bar{ A}^T\bar{ A}y+2\rho\Vert\bar{ A}yy^T\Vert^*\right\}\leq0.\label{Eq: optimization for Assumption positive semi-definite}
 			\end{eqnarray}
 			Now, we show that, given $\bar{ A}\in\mathbb{R}^n$, the optimal value of 
 			$
 			\sup_{y\in\mathbb{R}^n}\left\{-y^T\bar{ A}^T\bar{ A}y+2\rho\Vert\bar{ A}yy^T\Vert^*\right\}
 			$, denoted by $\tau^*$, is either $0$ or $+\infty$. 
 			\\
 			Clearly, the objective value of $y=0_{n}\in\mathbb{R}^n$ is $0$. Hence, $\tau^*\geq0.$ If for any $y\in\mathbb{R}^n$ the objective value is nonpositive, then $\tau^*=0$. Now, let us assume that there exists a $y\in\mathbb{R}^n$ such that the objective value is positive. Then for any $\alpha\in\mathbb{R}$, the objective value of $\alpha y$ is
 			$$
 			-\left(\alpha y\right)^T\bar{ A}^T\bar{ A}\left(\alpha y\right)+2\rho\Vert\bar{ A}\left(\alpha y\right)\left(\alpha y\right)^T\Vert^*=\alpha^2\left(-y^T\bar{ A}^T\bar{ A}y+2\rho\Vert\bar{ A}yy^T\Vert^*\right).
 			$$
 			Hence, the objective value of $\alpha y$ goes to $+\infty$ when $\alpha$ tends to $+\infty$. Therefore, the optimal value in this case is $+\infty.$}
 	\end{proof}
 	
 	\ahmadreza{Proposition \ref{Proposition: psd to sigma norm} provides an unconstrained optimization problem equivalent to checking Assumption \eqref{assumption: positive semidefinitness of a part}. 
 	}
 	\ahmadreza{In the next proposition, we show how one can use Proposition \ref{Proposition: psd to sigma norm} to check Assumption \eqref{assumption: positive semidefinitness of a part} for box uncertainty sets.}
 	\begin{proposition}\label{Prop:assumption psd for all with DC programming}
 		\ahmadreza{Let	$\uncertset=\left\{\Delta\in\mathbb{R}^{n\times n}:\ \Vert\Delta\Vert_\infty\leq\rho \right\}$. For this uncertainty set
 			Assumption \eqref{assumption: positive semidefinitness of a part}  holds if  and only if the following optimization probem has a nonnegative optimal value:
 			\begin{equation}\label{Eq:equavalent DC to assumption semi-definite for all}
 			\min_{y\in\mathbb{R}^n}\left\{y^T\bar{ A}^T\bar{ A}y-2\rho\Vert \bar{ A} y\Vert_1:\ \Vert y\Vert_1=1 \right\}.
 			\end{equation}}
 	\end{proposition}
 	\begin{proof}{Proof.} 
 		\ahmadreza{Using Proposition \ref{Proposition: psd to sigma norm}, Assumption \eqref{assumption: positive semidefinitness of a part}  holds if  and only if
 			\begin{equation}\label{Eq:equavalent nonlinear to assumption semi-definite for all}
 			\sup_{y\in\mathbb{R}^n}\left\{-y^T\bar{ A}^T\bar{ A}y+2\rho\Vert\bar{ A}yy^T\Vert_1\right\}=\sup_{y\in\mathbb{R}^n}\left\{-y^T\bar{ A}^T\bar{ A}y+2\rho\Vert\bar{ A}y\Vert_1\Vert y\Vert_1\right\}=0.
 			\end{equation}
 			As it is mentioned in the proof of Proposition \ref{Proposition: psd to sigma norm}, the optimal value of \eqref{Eq:equavalent nonlinear to assumption semi-definite for all} is either zero or $+\infty$.
 			Therefore, we know if Assumption \eqref{assumption: positive semidefinitness of a part} does not hold then the optimal value of \eqref{Eq:equavalent nonlinear to assumption semi-definite for all} is $+\infty$ and hence there exists $\bar{ y}\in\mathbb{R}^n$ such that 
 			$$\bar{ y}^T\bar{ A}^T\bar{ A}\bar{ y}-2\rho\Vert \bar{ A}\bar{ y}\Vert_1\Vert \bar{ y}\Vert_1 <0.$$
 			Let us define $\hat{y}:=\frac{\bar{ y}}{\Vert \bar{ y}\Vert_1}$. So, we have:
 			$$
 			\hat{ y}^T\bar{ A}^T\bar{ A}\hat{ y}-2\rho\Vert \bar{ A}\hat{ y}\Vert_1= \hat{ y}^T\bar{ A}^T\bar{ A}\hat{ y}-2\rho\Vert \bar{ A}\hat{ y}\Vert_1\Vert \hat{ y}\Vert_1=\frac{1}{\Vert \bar{ y}\Vert_1^2}\left(\bar{ y}^T\bar{ A}^T\bar{ A}\bar{ y}-2\rho\Vert \bar{ A}\bar{ y}\Vert_1\Vert \bar{ y}\Vert_1 \right)<0.
 			$$
 			So, in  \eqref{Eq:equavalent DC to assumption semi-definite for all} the objective value of the feasible solution $\hat{y}$ is negative. Hence, the optimal value of \eqref{Eq:equavalent DC to assumption semi-definite for all} is negative.
 			\\
 			Now, let us assume that Assumption \eqref{assumption: positive semidefinitness of a part} holds. By  Proposition \ref{Proposition: psd to sigma norm} for any $y\in \mathbb{R}^n$ we have that 
 			$-{ y}^T\bar{ A}^T\bar{ A}{ y}+2\rho\Vert \bar{ A}{ y}\Vert_1\Vert { y}\Vert_1 \leq 0,$ for any $y\in\mathbb{R}^n$.
 			Therefore, for any $y$ in the set $\mathcal{F}:=\left\{y\in\mathbb{R}^n:\hs\Vert { y}\Vert_1=1  \right\}$, we have ${ y}^T\bar{ A}^T\bar{ A}{ y}-2\rho\Vert \bar{ A}{ y}\Vert_1\geq0$ and hence \eqref{Eq:equavalent DC to assumption semi-definite for all} has a nonnegative optimal value.
 		}
 	\end{proof}

 	\ahmadreza{We emphasize that the optimization problem in \eqref{Eq:equavalent DC to assumption semi-definite for all} belongs to the class of  DC (Difference of Convex) optimization problems, for which an extensive literature exists (see, e.g., \cite{Ahmadi2018,Lipp2016,shen2016disciplined}).}
 	
 	\begin{remark}\label{Remark: approximation of assumption C}
 		\ahmadreza{ In this paper, we use \eqref{Eq:approximation of DC formulation of assumption psd for all}  to check Assumption \eqref{assumption: positive semidefinitness of a part}.  If the optimal solution of \eqref{Eq:approximation of DC formulation of assumption psd for all} is nonnegative, then the inequality in \eqref{Eq:equavalent DC to assumption semi-definite for all} and hence Assumption \eqref{assumption: positive semidefinitness of a part} hold. Furthermore, \eqref{Eq:approximation of DC formulation of assumption psd for all} can be seen as an  approximation of \eqref{Eq:equavalent DC to assumption semi-definite for all}.} \ahmadreza{The intuition behind this relaxation is that in \eqref{Eq:equavalent DC to assumption semi-definite for all} we want to minimize $\Eucledianorm{\bar{ A}y}$ and simultaneously maximizing $\Vert\bar{ A}y\Vert_1$. Due to symmetrical behaviour of norm functions, we only restrict the problem to be optimized on $\left\{y\in\mathbb{R}^n:\hs\bar{ A} y\geq0,\hs \Vert y\Vert_1=1\right\}$, and then relax it into  $\left\{y\in\mathbb{R}^n:\hs\bar{ A} y\geq0,\hs \Vert y\Vert_1\leq1\right\}$.} \hfill \qed
 	\end{remark}
 	\begin{remark}
 		\ahmadreza{	For a general uncertainty set, using the $\Omega$ from Assumption \eqref{assumption: upper bound on the uncertainty} and Proposition \ref{Prop:assumption psd for all, ellipsoid}, we know Assumption \eqref{assumption: positive semidefinitness of a part} holds if $\lambda_{\min}(\bar{ A}^T\bar{ A})\geq4\Omega^2$, where $\lambda_{\min}$ denotes the smallest eigenvalue of $\bar{ A}^T\bar{ A}.$
 			%
 			%
 		}\hfill \qed
 	\end{remark}

 	\section{\ahmadreza{A heuristic method to find a worst-case scenario in a norm approximation problem}}\label{Appendix:heuristic}
 	\ahmadreza{In this section, we describe a heuristic method to find a worst-case scenario in the uncertainty set $\mathcal{Z}$, defined in \eqref{Eq:uncertainty set with two constraints}, corresponding to  a given solution $y$. We assume, without loss of generality, that the components of $y$ are sorted such that their absolute values are descending.  }
 	
 	\ahmadreza{The idea behind the algorithm comes from the fact that for the box uncertainty set $$\left\{\Delta\in \mathbb{R}^{n\times n}: \;\Vert\Delta\Vert_\infty\leq \rho \right\}$$ and a solution $y$, the worst-case scenario $\Delta ^\square$ is defined by 
 		$$
 		\Delta^\square_{ij}=  sgn([Ay-b]_i) sgn(y_j)\rho,\hs i,j=1,...,n,
 		$$
 		where $sgn(.)$ is the sign function. To see this, we first rewrite $\Eucledianorm{(\bar{ A}+\Delta)y-b}$ as $\Eucledianorm{\Delta y+\bar{ A}y-b}$ and we recall that the worst-case scenario is a vertex of the box uncertainty set, which are matrices whose components are either $\rho$ or $-\rho$. We know that $\Delta$ maximizes  $\Eucledianorm{\Delta y+\bar{ A}y-b}$ if for any $i=1,...,n$, the value of $(\Delta y)_i$ is the highest and its sign is the same as the sign of $(\bar{ A}y-b)_i$. To make the value of $(\Delta y)_i$ the highest for any $i=1,...,n$, any component of $\Delta_{ij}$ should have the same sign as $y_j$, with $\vert \Delta_{ij}\vert=\rho$,  $j=1,...,n.$ Hence, $\Delta ^\square$ is the maximizer of $\Eucledianorm{(\bar{ A}+\Delta)y-b}$ over the box uncertainty set. }
 	
 	\ahmadreza{So, in the definition of $\uncertset$, if $B^1=B^2=\zero{n\times n}$, then $\Delta ^\square$ is the worst-case scenario. Also, if there exists $i\in\{1,...,n\}$ such that the $i$th column of $B^1$ and $B^2$ are zero, then for any $j=1,...,n$, the component in the $j$th row and $i$th column of the worst-case scenario is $\Delta ^\square_{ij}$. So, from now on we assume that $B^1$ and $B^2$ do not have any zero columns in common. }
 	
 	\ahmadreza{Let $\tilde{\Delta}=\zero{n\times n}$. As $y_1$ has the largest absolute value, changing the first column of $\tilde{\Delta}$ may result in a high  increase in the value of  $\Eucledianorm{(\bar{ A}+\tilde{\Delta})y-b}$. So, we start with $j=1$ and $i=1$. If $B^1_{ij}$ (or $B^2_{ij}$) is $1$, then we change $\tilde{\Delta}_{ij}$ to $sgn([\bar{ A}y-b]_i)_{n,n} sgn(y_j)\rho$. We increase $i$ by one and continue the same procedure. We also make sure that the number of changes happened because of $B^1_{ij}=1$ (or $B^2_{ij}=1$) does not exceed $K$. If we finish the procedure with $i=n$, then we increase $j$ by one and reset $i=1$.  This procedure guarantees that the resulting scenario $\tilde{\Delta}$ is feasible. In the next example, we illustrate how the algorithm works.
 		\begin{example}\label{Example:illustration of worst-case algorithm}
 			Let $n=3$, 
 			$$ \bar{ A}=\left[\begin{matrix}
 			2.8 &3.2 &5.1 \\-2.5& 3.6&0\\-1.5&2.7&3.0
 			\end{matrix}\right],\enskip B^1=\left[\begin{matrix}
 			1 &0&1 \\0&1&1\\1&1&1
 			\end{matrix}\right], \enskip B^2=\left[\begin{matrix}
 			1 &1&1 \\1&0&1\\1&0&1
 			\end{matrix}\right],\enskip b=\left[\begin{matrix}
 			2\\3\\1
 			\end{matrix}\right],\enskip y=\left[\begin{matrix}
 			-3\\2\\-1
 			\end{matrix}\right],
 			$$ 
 			$K=5,$ and $\rho=0.2.$ For this example, $\bar{ A}y-b=\left[\begin{matrix}-9.1&11.7&5.9\end{matrix}\right]^T$ and
 			$$
 			\Delta ^\square=\left[\begin{matrix}
 			0.2  & -0.2 & 0.2 \\
 			-0.2 & 0.2  & -0.2  \\
 			-0.2  & 0.2 & -0.2
 			\end{matrix}\right].	$$
 			Let $C_1$ and $C_2$ be the changes occurring because of $B^1$ and $B^2$, respectively.  For $\tilde{\Delta},$ we start with $i,j=1$. Since both $B^1_{11}=B^1_{11}=1$, so $\tilde{\Delta}_{11}=0.2$, $C_1=1$, and $C_2=1$. Since $K=5$, we have $C_1\leq K,$ and $C_2\leq K.$ So, we continue and increase $i$ to $2$. We have $B^2_{21}=1$ while $B^1_{21}=0$, so we change $\tilde{\Delta}_{21}=-0.2,$ and $C_2=2$. By continuing this procedure, we end up with 
 			$$
 			\tilde{\Delta}=\left[\begin{matrix}
 			0.2  & -0.2 & 0.2 \\
 			-0.2 & 0.2  & 0    \\
 			-0.2  & 0.2 & 0
 			\end{matrix}\right].
 			$$ 
 			For this example, the worst-case value of $\Eucledianorm{(A+\Delta)y-b}$ over $\uncertset$ is 17.78 (obtained using SCIP 5.1 \cite{gleixner2017scip}) while $\Eucledianorm{(A+\tilde{\Delta})y-b}=17.75$.  \hfill \qed
 		\end{example}
 		As one can also see from Example \ref{Example:illustration of worst-case algorithm}, the algorithm proposed in this section is a heuristic and the obtained scenario may not be the worst-case scenario.}
 	
 \end{appendices}
 
  \bibliographystyle{plain}
  \bibliography{Ref}


\end{document}